%% file: main.tex
\documentclass[12pt,a4paper]{amsart}
\usepackage[margin=30mm,includeheadfoot]{geometry}
\usepackage{graphicx} 
\usepackage{amsmath}
\usepackage{amssymb}
\usepackage{amscd}
\usepackage{amsthm}
\usepackage{needspace}
\usepackage{enumerate}
\usepackage{xcolor}
\usepackage[all]{xy}
\usepackage{url}
\usepackage[hidelinks,breaklinks]{hyperref}
\usepackage[lite,abbrev,alphabetic]{amsrefs}
\usepackage{ragged2e}
\AddToHook{env/biblist/begin}{\RaggedRight}
\usepackage{cleveref}
\usepackage{comment}
\usepackage{todonotes}
\usepackage{tikz-cd}
\usetikzlibrary{arrows.meta}
\usepackage{quiver}
\usepackage{enumitem}

\makeatletter
\renewcommand{\subsection}{\@startsection{subsection}{2}%
  \z@{.7\linespacing\@plus.5\linespacing}{.35\linespacing}%
  {\normalfont\bfseries\raggedright}}
\makeatother

\theoremstyle{definition}
\newtheorem{definition}{Definition}[section]

\theoremstyle{plain}
\newtheorem{theorem}[definition]{Theorem}
\newtheorem{proposition}[definition]{Proposition}
\newtheorem{lemma}[definition]{Lemma}
\newtheorem{corollary}[definition]{Corollary}

\newtheorem*{theorem*}{Theorem}

\newtheorem*{maintheorem}{Main theorem}

\theoremstyle{remark}
\newtheorem{remark}[definition]{Remark}
\newtheorem*{remark*}{Remark}

\crefname{theorem}{Theorem}{Theorems}
\crefname{proposition}{Proposition}{Propositions}
\crefname{lemma}{Lemma}{Lemmas}
\crefname{corollary}{Corollary}{Corollaries}
\crefname{conjecture}{Conjecture}{Conjectures}
\crefname{hypothesis}{Hypothesis}{Hypotheses}
\crefname{remark}{Remark}{Remarks}
\crefname{remark*}{Remark}{Remarks}
\crefname{maintheorem}{Main theorem}{Main theorems}
\crefname{condition}{Condition}{Conditions}
\crefname{example}{Example}{Examples}

\def\C{\mathbb{C}}
\def\Q{\mathbb{Q}}

\def\GL{\mathrm{GL}}
\def\PGL{\mathrm{PGL}}
\def\SL{\mathrm{SL}}

\def\SO{\mathrm{SO}}
\def\GSO{\mathrm{GSO}}

\def\Spin{\mathrm{Spin}}
\def\PGSO{\mathrm{PGSO}}
\def\PGSp{\mathrm{PGSp}}

\DeclareMathOperator{\diag}{diag} 

\DeclareMathOperator{\Hom}{Hom}
\DeclareMathOperator{\End}{End}
\DeclareMathOperator{\Aut}{Aut}

\DeclareMathOperator{\Cent}{Cent}
\DeclareMathOperator{\Irr}{Irr}
\DeclareMathOperator{\Std}{Std}
\DeclareMathOperator{\tr}{tr}
\DeclareMathOperator{\id}{id}
\DeclareMathOperator{\Sym}{Sym}

\newcommand{\cH}{\mathcal{H}}

\newcommand{\cS}{\mathcal{S}}
\newcommand{\cT}{\mathcal{T}}

\numberwithin{equation}{section}

\newcommand{\temp}{\mathrm{temp}}
\newcommand{\gen}{\mathrm{gen}}

\newcommand{\rG}{\mathrm{G}}
\newcommand{\Art}{\mathrm{Art}}
\newcommand{\ab}{\mathrm{ab}}

\newcommand{\abs}[1]{\lvert #1 \rvert}

\newcommand{\GSpin}{\mathrm{GSpin}}

\newcommand{\Trans}{\mathrm{Trans}}
\newcommand{\St}{\mathrm{St}}
\newcommand{\sn}{\mathrm{sn}}

\DefineSimpleKey{bib}{arxiv}
\newcommand\myurl[1]{\url{#1}}
\newcommand\arxiv[1]{available at \href{https://arxiv.org/abs/#1}{arXiv:#1}}

\BibSpec{article}{%
  +{}  {\PrintAuthors}                {author}
  +{,} { \textit}                     {title}
  +{.} { }                            {part}
  +{:} {\textit}                     {subtitle}
  +{,} {\PrintContributions}         {contribution}
  +{,} {\PrintConference}            {conference}
  +{}  {\PrintBook}                   {book}
  +{,} { }                            {booktitle}
  +{,} { }                            {journal}
  +{}  { \textbf}                     {volume}
  +{}  { \PrintDate}                  {date}
  +{,} { \issuetext}                  {number}
  +{,} { }                            {pages}
  +{,} { available at \eprint}        {eprint}
  +{;} { \arxiv}               {arxiv}
  +{.\!\!}  {\SentenceSpace \PrintReviews} {review}
  +{,} { \PrintDOI}                   {doi}
}

\allowdisplaybreaks[1]
\hypersetup{
  bookmarksdepth=2,bookmarksnumbered=true,bookmarksopen=false,
  pdftitle={Arthur packets for G2 over p-adic fields},
  pdfauthor={Yugo Takanashi},
  pdfsubject={2020 Mathematics Subject Classification: 22E50, 11F70, 20G41},
  pdfkeywords={Arthur packets, exceptional groups, endoscopy, triality, Aubert duality}
}

\title[Arthur packets for $\rG_2$ over $p$-adic fields]{Arthur packets for $\rG_2$ over $p$-adic fields}
\author{Yugo Takanashi}
\address{Yugo Takanashi\\Department of Mathematics, Graduate School of Science, Kyoto University, Kitashirakawa Oiwake-cho, Sakyo-ku, Kyoto 606-8502, Japan}
\email{takanashi.yugo.7s@kyoto-u.ac.jp}
\thanks{This work was supported by JSPS KAKENHI Grant Number JP26KJ0167.}
\date{}

\begin{document}
\begin{abstract}
      We construct Arthur packets for every non-tempered Arthur parameter of the split group $\rG_2$ over a $p$-adic field and prove their ordinary endoscopic character identities.
      We identify the corresponding distinguished Xu packets on $\PGSO_8$ and determine their triality-twisted character identities.
      The construction uses the twisted Aubert duality of Cheng--Kaletha.
      We also extend the construction and character identities to long-root and short-root parameters arising from generic unitary representations of $\PGL_2$.
\end{abstract}

\maketitle

\section*{Introduction}

Let $F$ be a finite extension of $\Q_p$.
Gan and Savin constructed the local Langlands correspondence for the split group $\rG_2$ over $F$ \cite{GanSavin2023-LocalLanglandsConjectureG2}.
The author proved the ordinary endoscopic character relations for its tempered packets and the triality-twisted character relations for their lifts to $\PGSO_8$ \cite{TakanashiEndoscopyLocal}.
Several families of local non-tempered Arthur packets were previously constructed using theta correspondences: the cubic unipotent packets arose from the work of Huang--Magaard--Savin \cite{HuangMagaardSavin1998-UnipotentRepresentationsG2} and Gan--Gurevich--Jiang \cite{GanGurevichJiang2002-CubicUnipotentArthurParameters}, and the short-root packets were constructed through the Rallis--Schiffmann correspondence by Gan--Gurevich \cite{GanGurevich2006-NontemperedAPacketsG2}.

In this paper, we construct packets for every non-tempered Arthur parameter of $\rG_2(F)$ and prove both their ordinary endoscopic character identities and the triality-twisted identities for their lifts to $\PGSO_8(F)$.
We explicitly identify the distinguished Xu packet attached to each lifted parameter and express the triality-twisted transfer of the stable virtual character on $\rG_2(F)$ in terms of the canonical twisted characters of that packet.
Together with the tempered results of \cite{TakanashiEndoscopyLocal}, this gives an endoscopic description of the local Arthur packets of $\rG_2(F)$.

\begin{maintheorem}[Informal; see Theorem~\ref{thm:main-ecr}]\label{thm:main-intro}
      For every non-tempered Arthur parameter $\psi$ of $G_2(F)$, we construct a finite packet $\Pi_\psi$ consisting of irreducible representations of $G_2(F)$ with a surjective map
      \begin{align*}
            \Pi_{\psi} \twoheadrightarrow \Irr_{\C}(\cS_{\psi}),
      \end{align*}
      and a stable virtual character $S\Theta_\psi$.
      These packets have the following properties:
      \begin{itemize}
            \item They satisfy the endoscopic character identities for the endoscopic groups of $G_2(F)$ and the triality-twisted endoscopic character identities for $\PGSO_8(F)$.

            \item The triality-twisted transfer of $S\Theta_\psi$ to $\PGSO_8$ is an explicit signed sum of canonical twisted characters in the corresponding distinguished Xu packet constructed by Xu in \cite{Xu2024-ArthurPacketsQuasisplitGSp2nGO2n}.
            This property characterizes the local Arthur packets of $G_2(F)$ uniquely.

            \item For non-tempered co-tempered parameters, an explicit sign formula relates their stable characters to the Aubert duals of the corresponding tempered stable characters.
            For co-tempered parameters, the map above is also bijective, and the compatibility with Aubert duality also gives a characterization of the packets.
      \end{itemize}
\end{maintheorem}

For co-tempered parameters, we apply Aubert duality to the tempered packet obtained by interchanging the Deligne and Arthur $\SL_2$ factors.
We use the twisted Aubert duality of Cheng--Kaletha \cite{ChengKaletha2026-AubertBernsteinInvolutionsDisconnectedGroups} for the triality-twisted identities.
The mixed cases also use Aubert duality to reduce the problem to long-root parameters and compute the packet explicitly by using known computations on composition series.
For $\rG_2$, these calculations are facilitated by Qin's computations of Aubert duality \cite{Qin2026-AubertZelevinskyInvolutionG2HeckeAlgebras}, together with the composition series of Mui\'c \cite{Muic1997-UnitaryDualG2}.
Our use of Aubert duality is inspired by Arthur's treatment \cite{Arthur2013-EndoscopicClassificationRepresentations}*{Section 7.1}, \cite{Hiraga2004-FunctorialityZelevinskiInvolutions} and the work of Atobe--Gan--Ichino--Kaletha--M\'inguez--Shin \cite{AtobeGanEtAl2025-LocalIntertwiningRelationsCoTemperedPackets}.

For archimedean fields, the short root Arthur packets are constructed by Gan--Gurevich \cite{GanGurevich2006-NontemperedAPacketsG2}*{Sections 10 and 11} and the subregular packets are described in Gan--Gurevich--Jiang \cite{GanGurevichJiang2002-CubicUnipotentArthurParameters}*{Section 6} via theta correspondences.
We will treat the archimedean case in a forthcoming paper, using an approach similar to that of this article.

The compatibility of our construction of Arthur packets with the Huang--Magaard--Savin, Gan--Gurevich--Jiang, and Gan--Gurevich constructions will be established in a forthcoming paper, as a step towards the endoscopic classification of discrete automorphic representations of $\rG_2$.
In recent work, Cunningham--Hazeltine--Liu--Lo--Ray--Xu--Sch\"urmann develop a framework to study the Kottwitz--Shelstad transfer for ABV-packets.
They also suggest applying their framework to the triality-twisted endoscopy relating $\rG_2$ and $\PGSO_8$; see \cite{CunninghamEtAl2026-GeometrizationKottwitzShelstadTransfer}*{Section 1.5}.
It may be interesting to compare our packets with the corresponding ABV-packets.

\subsection*{Outline of the paper}
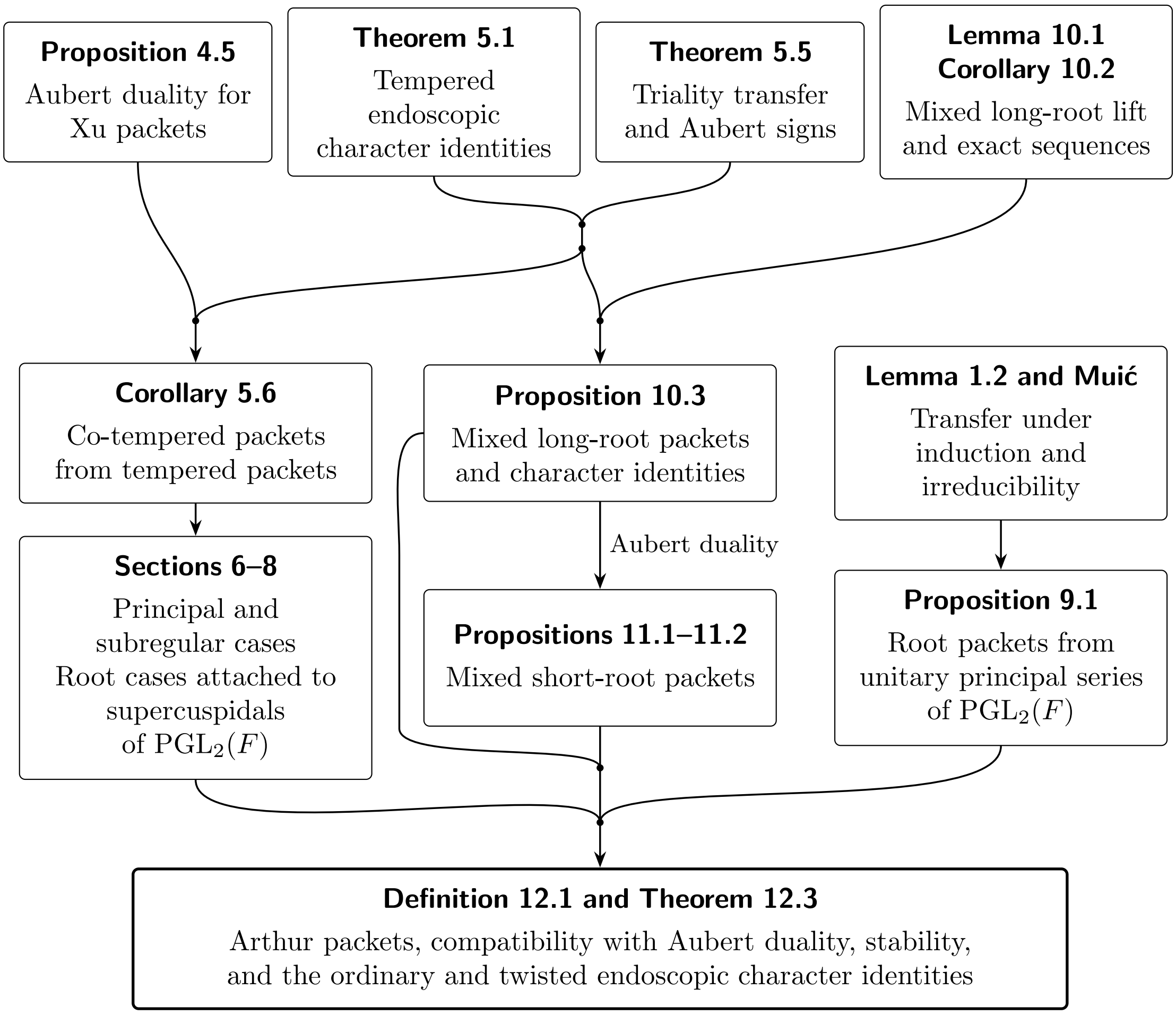
\begin{figure}[htbp]
      \centering
      \input{proof-of-main-theorem}
      \caption{Proof of Theorem~\ref*{thm:main-ecr}. Arrows point from inputs to conclusions. Joining curves indicate inputs used together.}
      \label{fig:proof-of-main-theorem}
\end{figure}

Figure~\ref{fig:proof-of-main-theorem} summarizes the proof of Theorem~\ref{thm:main-ecr}.
\begin{itemize}
      \item Section~\ref{sec:generalities} collects the general representation-theoretic notation, parameter conventions, transfer normalizations, and general results on Aubert duality.
      \item Sections~\ref{sec:preliminaries-G2-PGSO8}--\ref{sec:twisted-endoscopy} specialize these inputs to $\rG_2$ and the orthogonal groups, including the comparison with Qin's notation, the tempered packet identities, and the signs for triality transfer.
      \item Sections~\ref{sec:principal-parameter}--\ref{sec:mixed-short} compute the elements of A-packets and the character identities for each family of parameters.
      \item Section~\ref{sec:summary} states and proves the main theorem from these calculations.
      \item Section~\ref{sec:A-parameter-unbounded-ecr} also treats root parameters arising from unitary generic representations of $\PGL_2(F)$, including complementary series.
\end{itemize}

\Needspace{10\baselineskip}
\begin{remark*}
      Throughout this paper, we assume the twisted weighted fundamental lemma for the unramified twisted endoscopic Levi data that occur for $\PGSO_8$ with triality, as in \cite{TakanashiEndoscopyLocal}.
      This assumption enters through the tempered endoscopic character identities recalled in Theorem~\ref{thm:tempered-input}.
      For the twisted weighted fundamental lemma and its reduction, see \cite{MoeglinWaldspurger2016-StabilisationFormuleTracesTordue1}*{II.4.4} and \cite{AtobeGanEtAl2025-LocalIntertwiningRelationsCoTemperedPackets}*{Section 0.4}.
      In sufficiently large residue characteristic, the twisted weighted fundamental lemma follows from the Lie algebra and nonstandard weighted fundamental lemmas of Halleck-Dub\'e \cite{HalleckDube2026-WeightedFundamentalLemmaNonSplitGroups}*{Theorems 16.6 and 17.9}, as explained in Section~1.14 of the same paper.
\end{remark*}

\section{Generalities}\label{sec:generalities}

\subsection{Basic conventions}\label{subsec:basic-conventions}

Throughout, $F$ is a finite extension of $\Q_p$.
We fix a nontrivial additive character $\psi_F \colon F \to \C^\times$.
In this section, let $G$ denote a connected reductive group over $F$ and $P = MN$ a parabolic subgroup of $G$ with Levi factor $M$ and unipotent radical $N$.
This section establishes general conventions and tools used below.

We use Artin reciprocity $\Art_F \colon F^{\times} \to W_F^{\ab}$, which sends a uniformizer to an arithmetic Frobenius element.
A character $\chi$ of $F^{\times}$ is viewed as a character of $W_F$ via $\chi \circ \Art_F^{-1}$.
In particular, $\nu(w) = \abs{w}_F$.
We set $\nu_n=|\det|_F$ on $\GL_n(F)$.
For the group $\GSO_{n}$, let $\mathrm{sim} \colon \GSO_n(F) \to F^{\times}$ denote the similitude character.

\subsection{Representations and normalized inductions}\label{subsec:general-representations}

For a connected reductive $F$-group $G$, write $\Irr(G(F))$ for the isomorphism classes of irreducible admissible complex representations, and $\operatorname{JH}$ for irreducible constituents of a finite-length representation.
Let $R(G)$ be the Grothendieck group of finite-length smooth complex representations of $G(F)$, and write $\Theta_\pi(f)=\tr\pi(f)$ for $f\in C_c^\infty(G(F))$.
If a representation $\pi$ has a central character, let $\omega_\pi$ denote it.

Fix a minimal $F$-parabolic $P_{0, G}=M_{0, G}N_{0, G}$.
Standard parabolics contain $P_{0, G}$, and their standard Levi factors contain $M_{0, G}$.
We sometimes omit the subscript $G$ when the group is clear from context.
For a parabolic $P=MN$ of $G$, let $\delta_P$ be its modulus character and put
\[
      I_P^G(\sigma)=
      \operatorname{Ind}_{P(F)}^{G(F)}(\delta_P^{1/2}\sigma),
      \qquad
      r_P^G(\pi)=\delta_P^{-1/2}\otimes\pi_{N(F)}.
\]
Write $\overline P=M\overline N$ for the parabolic opposite to $P$.
Here $\sigma$ is extended trivially across $N(F)$, and $\pi_{N(F)}$ denotes the unnormalized Jacquet module.
When a standard Levi is displayed in an induction or quotient subscript, it denotes the standard parabolic with that Levi.
For orthogonal groups, $\tau\rtimes\rho$ denotes normalized induction from the displayed Levi $\GL_n\times\SO_m$, with $\tau$ on $\GL_n(F)$ and $\rho$ on $\SO_m(F)$.

Let $X_F^*(M)$ be the group of $F$-rational algebraic characters of $M$, and put
\[
      \mathfrak a_M=\Hom_{\mathbf Z}(X_F^*(M),\mathbf R).
\]
Write $\mathfrak a_M^*$ for its real dual and $\mathfrak a_{M,\C}^*=\mathfrak a_M^*\otimes_{\mathbf R}\C$.
Define $H_M:M(F)\to\mathfrak a_M$ by $\langle\xi,H_M(m)\rangle=\log|\xi(m)|_F$.
For $\lambda\in\mathfrak a_{M,\C}^*$, set
\[
      \nu_M^\lambda(m)=\exp\langle\lambda,H_M(m)\rangle.
\]
We also write $\nu^\lambda$ when $M$ is understood.
Thus $\nu_M^{s\xi}=|\xi|_F^s$ for $\xi\in X_F^*(M)$.
If $\sigma$ is irreducible tempered and $\lambda\in\mathfrak a_M^*$ is strictly positive relative to a standard $P$, then $I_P^G(\sigma\otimes\nu_M^\lambda)$ is a Langlands standard module.
Its unique irreducible quotient is denoted by $J_P^G(\sigma\otimes\nu_M^\lambda)$.
The notation $J_P^G(\eta)$ is used only when the displayed induction is such a standard module.
We also allow $P=G$ and $\lambda=0$, giving an irreducible tempered standard module.

We now specialize to $G = \GL_2$.
Let $B_2$ be the upper triangular Borel of $\GL_2$.
For arbitrary smooth characters $a,b$ of $F^\times$, put
\[
      \pi(a,b)=I_{B_2}^{\GL_2}(a\otimes b).
\]
Here $(a\otimes b)(\diag(x,y))=a(x)b(y)$.
This is the entire normalized principal series, which may be reducible; if $a,b$ are unitary, it is irreducible tempered.
Write $\St_G$ for the Steinberg representation of $G(F)$, $\St_n=\St_{\GL_n}$, and, for a smooth character $\chi$ of $F^\times$, $\St_2(\chi)=\St_2\otimes(\chi\circ\det)$, so $\St_2=\St_2(1)$.
A twist $\tau\nu_2^s$ means $\tau\otimes|\det|_F^s$.

We shall also use the standard $\GL_2$ exact sequence
\begin{equation}\label{eq:GL2-length-two}
      0\longrightarrow\St_2(\chi)\nu_2^{1/2}
      \longrightarrow\pi(\chi\nu,\chi)
      \longrightarrow(\chi\circ\det)\nu_2^{1/2}\longrightarrow0,
\end{equation}
valid for every smooth character $\chi$ of $F^\times$; see \cite{GanSavin2023-HoweDualityDichotomyExceptionalThetaCorrespondences}*{Section 3.1}.

\subsection{Parameters and infinitesimal characters}\label{subsec:parameters-and-infinitesimal-characters}

The parameter conventions needed here concern split connected reductive groups $G$.
We fix a pinned dual group $\widehat G$, write ${}^LG=\widehat G\times W_F$, and display only the $\widehat G$-component of each parameter; its $W_F$-component is the projection onto $W_F$.
Write $L_F=W_F\times\SL_2^D(\C)$, where $D$ denotes the Deligne factor.
A Langlands parameter is an admissible homomorphism $\varphi:L_F\to\widehat G$, algebraic on $\SL_2^D$; it is tempered when its restriction to $W_F$ is bounded.
Parameters are considered up to $\widehat G$-conjugacy.
Put $d_w=\diag(|w|_F^{1/2},|w|_F^{-1/2})$.

\begin{definition}\label{def:parameters}
      An Arthur parameter for $G$ is an admissible homomorphism
      \[
            \psi:W_F\times\SL_2^D(\C)\times\SL_2^A(\C)\longrightarrow\widehat G
      \]
      that is bounded on $W_F$ and algebraic on the two $\SL_2$ factors; $A$ denotes the Arthur factor.
      Its associated Langlands parameter and its swapped parameter are
      \begin{align}
            \phi_\psi(w,x)&=\psi(w,x,d_w),\label{eq:associated}\\
            \widehat\psi(w,x,y)&=\psi(w,y,x).\notag
      \end{align}
      We call $\psi$ tempered if it is trivial on $\SL_2^A$, and co-tempered if it is trivial on $\SL_2^D$; equivalently, $\widehat\psi$ is tempered.

      We also set
      \begin{align*}
            s_{\psi}=\psi(1, I_2, -I_2).
      \end{align*}
\end{definition}

Let $S_a=\Sym^{a-1}(\C^2)$.
In any fixed linear representation of $\widehat G$, a summand $\rho\boxtimes S_a^D\boxtimes S_b^A$ of $\psi$ contributes
\begin{equation}\label{eq:parameter-block}
      \bigoplus_{j=0}^{b-1}\rho\nu^{(b-1)/2-j}\boxtimes S_a^D
\end{equation}
to $\phi_\psi$.
Indeed, the diagonal torus acts on $S_b$ with weights $b-1,b-3,\ldots,1-b$.

The infinitesimal character of a Langlands parameter $\varphi$ is the $\widehat G$-conjugacy class of the homomorphism
\[
      \lambda_\varphi(w)=\varphi(w,d_w).
\]
For an Arthur parameter this gives
\begin{equation}\label{eq:infinitesimal-Arthur}
      \lambda_\psi(w)=\lambda_{\phi_\psi}(w)=\psi(w,d_w,d_w).
\end{equation}
In particular, swapping the two $\SL_2$ factors preserves the infinitesimal character.

The notion of infinitesimal character is compatible with parabolic induction and Jacquet module functors for quasi-split classical groups or similitude groups, see \cite{Xu2024-ArthurPacketsQuasisplitGSp2nGO2n}*{Corollary 5.5}.

\subsection{Endoscopic transfer, induction, and stability}\label{subsec:transfer-normalization}

Let $G$ be quasi-split, fix a pinning and its Whittaker datum using $\psi_F$, and let $a$ be a finite-order $F$-automorphism preserving the pinning.
Let $(H, s, \cH, \xi)$ be an endoscopic group of $(G, a)$.
In this paper, we only consider the case $$\cH=\widehat H\rtimes W_F$$ and $\xi$ is the natural embedding $\cH\hookrightarrow{}^LG$.

We distinguish the coset $G(F) \rtimes a$ from the disconnected group $G(F)\rtimes\langle a \rangle$.
For transfer on a fixed coset $G(F) \rtimes a$, we use the corrected Whittaker-normalized Artin factor $\Delta_A$ of \cite{AtobeGanEtAl2025-LocalIntertwiningRelationsCoTemperedPackets}*{Appendix F.7}, denoted $\Delta'_\lambda$ in \cite{KottwitzShelstad2012-SplittingInvariantsSignConventions}*{equation (5.5.2)}.
We use the discriminant-normalized orbital integrals and compatible measures of \cite{MoeglinWaldspurger2016-StabilisationFormuleTracesTordue1}*{I.2.4}.

Here $S^H(\delta,f^H)$ denotes a stable orbital integral and $I^{G \rtimes a}(\gamma,f)$ a twisted orbital integral.
Conjugacy on $G(F) \rtimes a$ is under $G(F)$, as in \cite{Kaletha2022-LocalLanglandsConjecturesDisconnectedGroups}*{Lemma 4.6}.
Put
\[
      \Delta(\delta,\gamma)
            =\Delta_A(\delta,\gamma)/\Delta_{\mathrm{IV}}(\delta,\gamma).
\]
Following \cite{TakanashiEndoscopyLocal}*{Section 1.3}, matching functions satisfy, for strongly $G$-regular $\delta$,
\[
      S^H(\delta,f^H)=d(a)^{1/2}
      \sum_{\gamma\leftrightarrow\delta}
      \frac{\Delta(\delta,\gamma)}{[Z_G(\gamma)(F):G_\gamma(F)]}
      I^{G \rtimes a}(\gamma,f).
\]
The sum is over rational conjugacy classes matching $\delta$, and $G_\gamma=Z_G(\gamma)^\circ$ is the connected centralizer used in $I^{G \rtimes a}$.
For the Lie algebra $\mathfrak t_G$ of the pinned maximal torus $T_G$, put
\[
      d(a)=
      \left|\det(1-a\mid\mathfrak t_G/\mathfrak t_G^a)\right|_F.
\]
Thus $d(1)=1$; for pinned triality of type $D_4$, $d(a)=|3|_F$.

Write $\operatorname{Tran}_H^{G \rtimes a}$ for spectral transfer from stable distributions on an endoscopic group $H$ to invariant distributions on $G(F) \rtimes a$.
Thus
\[
      \operatorname{Tran}_H^{G \rtimes a}(S)=T
      \quad\text{means}\quad T(f)=S(f^H).
\]
For ordinary endoscopy take $a=1$ and omit it from the notation.
A stable distribution is one that factors through stable orbital integrals.
In particular, stability of a character is tested on stable conjugacy classes of regular semisimple elements.

We specify the Levi data used for induction.
Let $\eta:{}^LH\hookrightarrow{}^LG$ be the fixed endoscopic embedding, and write $\dot s\in\widehat G\rtimes\widehat a$ for its endoscopic element.
Assume the ambient datum is elliptic and use embedded data with $\eta(\widehat H)=\Cent(\dot s,\widehat G)^\circ$ and $\eta({}^LH)\subset\Cent(\dot s,{}^LG)$, as in the applications.
For an $F$-parabolic $P_H=M_HN_H$, put
\[
      \widehat A_{M_H}=(Z(\widehat M_H)^{W_F})^\circ,
      \qquad \mathcal M=\Cent(\eta(\widehat A_{M_H}),{}^LG).
\]
The Weil action here is supplied by ${}^LM_H$.
By \cite{MoeglinWaldspurger2016-StabilisationFormuleTracesTordue1}*{I.3.4}, after conjugation $\mathcal M={}^LM$ for an $a$-stable Levi $M$ of $G$.
Choose $\lambda\in X_*(\widehat A_{M_H})$ in the open cone defining $\widehat P_H$, avoiding the kernels of the nonzero restrictions of $\widehat G$-roots to $\eta(\widehat A_{M_H})$.
Define $\widehat P$ by the nonnegative weights of $\eta\circ\lambda$, and let $P=MN$ be the corresponding parabolic.
Restriction gives
\[
      \eta_M:{}^LM_H\hookrightarrow{}^LM,
      \qquad
      \eta_M(\widehat M_H)=\Cent(\dot s_M,\widehat M)^\circ,
      \qquad \dot s_M\in\widehat M\rtimes\widehat a.
\]
Here $\dot s_M$ is the conjugate of $\dot s$ in the chosen Levi.
We assume this inherited datum is relevant, so matching strongly regular $F$-points exist.
Its auxiliary data, transfer factor, and measures are inherited from the ambient datum as in loc.~cit.
We call this the \emph{inherited Levi datum} and denote its spectral transfer by $\operatorname{Tran}_{M_H,\eta_M}^{M \rtimes a}$.

\begin{lemma}\label{lem:induction-transfer}
      For the inherited Levi datum above and a stable virtual character $S$ of $M_H(F)$, one has
      \[
            \operatorname{Tran}_{H,\eta}^{G \rtimes a}(I_{P_H}^H S)
              = I_P^G\bigl(\operatorname{Tran}_{M_H,\eta_M}^{M \rtimes a}(S)\bigr).
      \]
      Induction of distributions is defined by normalized constant terms; on the right it is evaluated on $G(F)\rtimes a$.
      For ordinary induction, $I_P^G S_M$ is stable whenever $S_M$ is a stable virtual character of $M(F)$.
\end{lemma}

\begin{proof}
      If $f^H$ is a transfer of $f$, write $f_P$ and $(f^H)_{P_H}$ for normalized constant terms.
      Descent with the inherited factors and measures gives
      \[
            [(f^H)_{P_H}]=[(f_P)^{M_H}]
      \]
      in the quotient of test functions on $M_H(F)$ by the common kernel of stable orbital integrals.
      This is \cite{MoeglinWaldspurger2016-StabilisationFormuleTracesTordue1}*{ I.4.11, condition (2)}, with the data of I.3.4.
      Evaluation on $S$ and the constant-term character formula prove the identity.
      For the last assertion, normalized constant terms induce a map from the stable quotient for $G$ to that for $M$; see loc.~cit., I.3.1.
      Its dual sends $S_M$ to the stable character $I_P^G S_M$.
\end{proof}

For the standard nonelliptic Levi datum we also use the following ordinary transfer identity.
Let $G$ be split, let $P=MN$ be standard, and take $H=M$ with the standard inclusion ${}^LM\hookrightarrow{}^LG$.
Choose $\dot s\in Z(\widehat M)^\circ$ with $\Cent(\dot s,\widehat G)^\circ=\widehat M$, and inherit the Whittaker normalization, auxiliary data, and measures from $G$.
Then, for every stable virtual character $S$ of $M(F)$,
\[
      \operatorname{Tran}_{M}^{G}(S)=I_P^G(S).
\]
This is the identity in \cite{Hiraga2004-FunctorialityZelevinskiInvolutions}*{Section 6, Lemma 6.1, p.~1642}, with the datum on $M$ equal to the identity datum.
Equivalently, normalized constant-term descent followed by evaluation on $S$ gives normalized induction.
This separate statement applies to the nonelliptic Levi datum; the ambient ellipticity assumption preceding Lemma~\ref{lem:induction-transfer} is not used.

For later use, we record a general criterion for automatic stability of Levi characters.

\begin{lemma}\label{lem:central-quotient-stability}
      Let $L$ be a product of general linear groups over $F$, and let $Z\subset L$ be a central split torus.
      Every finite-length character of $(L/Z)(F)$ is stable.
\end{lemma}

\begin{proof}
      If $T\subset L/Z$ is a maximal torus, its inverse image $T_L$ in $L$ is a product of multiplicative groups of finite etale $F$-algebras.
      Hence $H^1(F,T_L)=H^1(F,L)=1$.
      The central-extension cohomology sequences show that
      \[
            \ker\bigl(H^1(F,T)\longrightarrow H^1(F,L/Z)\bigr)=1,
      \]
      and the connecting map from $H^1(F,T)$ to $H^2(F,Z)$ is injective, and is compatible with the connecting map for $L/Z$.
      Thus a strongly regular semisimple stable conjugacy class in $L/Z$ contains only one rational conjugacy class.
      Every invariant character is therefore stable.
\end{proof}

\subsection{Canonical extensions for an $S_3$-action}\label{subsec:canonical-extensions-s3-action}

Let $G$ be a connected reductive group over $F$ with an action of $S_3$ by $F$-automorphisms.
Write $A_3 \subset S_3$ for the alternating subgroup and choose a generator $a$ of $A_3$.

\begin{proposition}\label{prop:canonical}
      If $\pi$ is an admissible representation of $G(F)$ such that $\End_{G}(\pi)=\mathbb{C}$ and $\pi \cong \pi^g$ for $g \in S_3$, then there exists an extension $\pi^{+}$ of $\pi$ to $G(F) \rtimes S_3$.
      The restriction of $\pi^{+}$ to $G(F) \rtimes A_3$ is uniquely determined by $\pi$ up to isomorphism.
\end{proposition}

\begin{proof}
      Choose operators $U, V$ implementing generators of orders $3$ and $2$ of $S_3$.
      We may assume that $U^3 = 1$ and $V^2 = 1$.
      Schur's lemma implies that
      \begin{align*}
            VUV^{-1} = c U^{-1}
      \end{align*}
      for some $c \in \mathbb{C}^{\times}$.
      Then, we have $c^3 = 1$.
      By replacing $U$ by $cU$, we have $VUV^{-1} = U^{-1}$.
      Thus, we obtain a representation of $G(F) \rtimes S_3$ on the same space as $\pi$.
      Both the trivial and sign characters of $S_3$ are trivial on $A_3$, so the restriction of $\pi^{+}$ to $G(F) \rtimes A_3$ is uniquely determined by $\pi$ up to isomorphism.
\end{proof}

For a finite-length representation $V$ of $G(F) \rtimes S_3$, we use the twisted character on $G(F) \times a$ defined by
\begin{align*}
      \Theta_V(f) = \tr \left( \int_{G(F)} V(g \rtimes a) f(g) dg \right)
\end{align*}
for $f \in C_c^{\infty}(G(F))$.

\begin{corollary}
      If $\pi$ is an admissible representation of $G(F) \rtimes S_3$ of length $2$ as a $G(F)$ representation with an irreducible socle $S$ and cosocle $C$ as $G(F)$ representations, then we have
      \begin{align*}
            \Theta_{\pi}(f)|_{G(F) \times a} = \Theta_{\widetilde{S}}(f) + \Theta_{\widetilde{C}}(f)
      \end{align*}
\end{corollary}

\subsection{Aubert involution and endoscopic transfer}\label{subsec:aubert-involution-and-endoscopic-transfer}

For a connected reductive $F$-group $G$ put
\begin{equation}\label{eq:Aubert-operator}
      \mathcal D_G=\sum_{P=MN\supset P_{0,G}}
         (-1)^{\dim(A_M/A_G)}I_P^G r_P^G.
\end{equation}
Here $A_M$ is the maximal split torus in the center of $M$.
For irreducible $\pi$ with supercuspidal support $(L,\sigma)$, put $\epsilon_G(\pi)=(-1)^{\dim(A_L/A_G)}$.

\begin{proposition}\label{prop:Aubert}
      There is a positive irreducible Aubert dual $\widehat\pi$ such that
      \[
            \mathcal D_G[\pi]=\epsilon_G(\pi)[\widehat\pi].
      \]
      On the subcategory whose irreducible subquotients have supercuspidal support on Levis of a fixed relative semisimple rank, the positive dual extends to an exact covariant functor.
      It is involutive on isomorphism classes and satisfies
      \[
            \widehat{I_P^G(\sigma)}\simeq I_{\overline P}^G(\widehat\sigma).
      \]
      For supercuspidal $\pi$, one has $\widehat\pi=\pi$.
\end{proposition}

\begin{proof}
      Let $r$ be the relative semisimple rank of $G$ and let $t$ be that of a Levi in the supercuspidal support.
      Give the term indexed by a subset $J$ of the relative simple roots degree $|J|$ in the Aubert complex $C_\pi$.
      Since the exponent in \eqref{eq:Aubert-operator} is $r-|J|$, one has $\mathcal D_G[\pi]=(-1)^r\chi(C_\pi)$.
      The cohomology is concentrated in the support degree $t$ and equals the positive dual.
      Hence
      \[
            \mathcal D_G[\pi]=(-1)^{r+t}[\widehat\pi]
                 =(-1)^{r-t}[\widehat\pi]
                 =\epsilon_G(\pi)[\widehat\pi].
      \]
      We use Aubert's irreducibility and involutivity theorems and the compatibility with induction in the formulation recalled in \cite{AtobeGanEtAl2025-LocalIntertwiningRelationsCoTemperedPackets}*{ Theorems B.2.1 and B.2.3}.
      For exactness, see also the homological duality of Schneider--Stuhler \cite{SchneiderStuhler1997-RepresentationTheorySheaves}*{Section III.3}.
      For a supercuspidal representation every proper Jacquet module is zero, so only $P=G$ remains in \eqref{eq:Aubert-operator} and the sign is $+1$.
\end{proof}

We write $D_A$ for this positive Aubert functor, on the subcategories of fixed relative support-Levi rank in Proposition~\ref{prop:Aubert}.

For a split semisimple group $G$ and a finite-order pinning-preserving automorphism $a$, we also use the twisted operator
\begin{equation}\label{eq:twisted-Aubert-operator}
      \mathcal D_G^a=
       \sum_{\substack{P=MN\supset P_{0,G}\\a(P)=P}}
         (-1)^{\dim A_M^a}I_P^G r_P^G.
\end{equation}
Here $P_{0,G}$ is the pinned Borel.
The exponent is the dimension of the identity component of the $a$-fixed subgroup of $A_M$, and the induced operators are evaluated on $G(F)\rtimes a$.
Since $G$ is semisimple, $\mathcal D_G^1=\mathcal D_G$.

Compatibility with ordinary endoscopic transfer is due to Hiraga \cite{Hiraga2004-FunctorialityZelevinskiInvolutions}*{Theorems 6.6--6.7}.
In our convention, for a split endoscopic group $H$ of a split $G$ it reads
\[
      \mathcal D_G\operatorname{Tran}_H^G
        =(-1)^{\dim A_H-\dim A_G}
            \operatorname{Tran}_H^G\mathcal D_H.
\]
Indeed, Hiraga's operator is $(-1)^{r(G)}\mathcal D_G$, where $r(G)$ is the semisimple split rank, and his transfer sign is $(-1)^{a(G)-a(H)}$, where $a(G)=r(G)+\dim A_G$.
In particular, the displayed sign is $+1$ for split semisimple groups.
We use Xu's extension to twisted endoscopy \cite{Xu2017-MoeglinParametrizationArthurPackets}*{Appendix A}.

\begin{proposition}\label{prop:Xu-transfer}
      Let $H$ be an endoscopic group of $(G, a)$, where the groups are split semisimple and $a$ preserves a pinning.
      Assume the datum is realized by an $L$-embedding as in \cite{Xu2017-MoeglinParametrizationArthurPackets}*{Appendix A}, with compatible transfer normalizations.
      Then, on stable virtual characters,
      \begin{equation}\label{eq:Xu-transfer}
            \mathcal D_G^a\operatorname{Tran}_H^{G \rtimes a}
              =\operatorname{Tran}_H^{G \rtimes a}\mathcal D_H.
      \end{equation}
      The Aubert operator on the right preserves stable virtual characters.
\end{proposition}

\begin{proof}
      We give the indexing in Xu's Jacquet-transfer identity.
      Conjugate the embedding $\eta$ so that $\eta(\widehat T_H)=(\widehat T_G^{\widehat a})^\circ$ and $\eta(\widehat B_H)\subset\widehat B_G$, and identify ${}^LH$ with its image.
      For an $a$-stable standard $P=MN$, put
      \[
            W_{G, a}=W(\widehat G,\widehat T_G)^{\widehat a},
            \quad W_{M,a}=W(\widehat M,\widehat T_G)^{\widehat a},
            \quad W_H=W(\widehat H,\widehat T_H)\subset W_{G, a}.
      \]
      Choose a torus $A\subset(\widehat T_G^{\widehat a})^\circ$ with ${}^LM=\Cent(A,{}^LG)$, and define
      \[
            \begin{aligned}
                  \Omega_H(M)&=\{w\in W_{G, a}:
                        \Cent(w(A),{}^LH)\longrightarrow W_F\text{ is surjective}\},\\
                  \mathcal W_H(M)&=W_H\backslash\Omega_H(M)/W_{M,a}.
            \end{aligned}
      \]
      For $w$ in this finite set, represented by $g\in N_{\widehat G}(\widehat T_G)$, let $P_{H,w}=M_{H,w}N_{H,w}$ be standard in $H$ with dual parabolic and Levi $\widehat H$-conjugate to
      \[
            g\,{}^LP\,g^{-1}\cap{}^LH,
            \qquad g\,{}^LM\,g^{-1}\cap{}^LH.
      \]
      Choose $h_w\in\widehat H$ carrying the standard dual parabolic and Levi to these two intersections.
      In the identification of ${}^LH$ with its image, the inherited embedding is
      \[
            \eta_{M,w}=\operatorname{Int}(g^{-1}h_w)|_{{}^LM_{H,w}}:
                              {}^LM_{H,w}\hookrightarrow{}^LM,
      \]
      up to $\widehat M$-conjugacy.
      With inherited normalizations, Xu's identity is
      \[
            r_P^G\operatorname{Tran}_{H,\eta}^{G \rtimes a}(S)
             =\sum_{w\in\mathcal W_H(M)}
                \operatorname{Tran}_{M_{H,w},\eta_{M,w}}^{M\rtimes a}
                                  (r_{P_{H,w}}^H S).
      \]
      Each double coset is counted once.
      Induction-transfer holds separately for each inherited embedding, by \cite{Xu2017-MoeglinParametrizationArthurPackets}*{ Appendix A, diagram (A.2) and the identities immediately following it}.
      This termwise identity also covers the inherited data that are nonelliptic in $M$; it is the statement in Xu's cited formula.
      Thus
      \[
            I_P^G r_P^G\operatorname{Tran}_{H,\eta}^{G \rtimes a}(S)
             =\operatorname{Tran}_{H,\eta}^{G \rtimes a}
                  \left(\sum_{w\in\mathcal W_H(M)}
                                I_{P_{H,w}}^H r_{P_{H,w}}^H S\right).
      \]
      Multiply by $(-1)^{\dim A_M^a}$ and sum over $P$.
      For each standard $P_H=M_HN_H$, let $a_{M_H,H,M}$ count the double cosets $w\in\mathcal W_H(M)$ in the preceding sum for which $P_{H,w}=P_H$.
      The coefficient on the right is
      \[
            \sum_P(-1)^{\dim A_M^a}a_{M_H,H,M}
                  =(-1)^{\dim A_{M_H}},
      \]
      by \cite{Xu2017-MoeglinParametrizationArthurPackets}*{Proposition A.1, equation (A.3)}.
      This is exactly $\mathcal D_H$ because $H$ is semisimple.
      This proves \cite{Xu2017-MoeglinParametrizationArthurPackets}*{diagram (A.1)}, whose left vertical arrow is the operator on stable distributions.
\end{proof}

\section{Preliminaries for $\rG_2$ and $\PGSO_8$}\label{sec:preliminaries-G2-PGSO8}

We now fix the root conventions and collect the facts needed for the groups $\rG_2$ and $\PGSO_8$.

\subsection{Explicit description}\label{subsec:pinned-groups}

Throughout the rest of the paper, $\rG_2$ denotes the split exceptional group, and we put
\[
      G=\SO_8,\qquad G_{\mathrm{ad}}=\PGSO_8=G/Z(G),\qquad G_{\mathrm{sim}}=\GSO_8.
\]
All these groups are split over $F$.
The notation $H$ is reserved for an endoscopic group of $\rG_2$.

We fix compatible pinnings and Whittaker data, with Haar measures and transfer-factor normalizations as in Section~\ref{subsec:transfer-normalization}.
We use the standard hyperbolic basis $e_1, \dots, e_8$ for the isotropic quadratic form on $F^8$ defining $G$ and $G_{\mathrm{sim}}$ to describe the standard maximal torus and Borel.

Write
\[
      \begin{aligned}
            B_{\rG_2}&=T_{\rG_2}U_{\rG_2},& B_G&=T_GU_G,\\
            B_{G_{\mathrm{ad}}}&=T_{G_{\mathrm{ad}}}U_{G_{\mathrm{ad}}},& B_{G_{\mathrm{sim}}}&=T_{G_{\mathrm{sim}}}U_{G_{\mathrm{sim}}}
      \end{aligned}
\]
for the pinned Borel subgroups; these are the minimal parabolics $P_{0,\rG_2}$, $P_{0,G}$, $P_{0,G_{\mathrm{ad}}}$, and $P_{0,G_{\mathrm{sim}}}$.
Let
\[
      j:G\hookrightarrow G_{\mathrm{sim}},\qquad
      p:G_{\mathrm{sim}} \longrightarrow G_{\mathrm{ad}},\qquad q=p\circ j:G\longrightarrow G_{\mathrm{ad}}
\]
be the inclusion and the quotients by the scalar centers.
We use $q^*$ for pullback of representations of $G_{\mathrm{ad}}(F)$ to $G(F)$ and $p^*$ for inflation to $G_{\mathrm{sim}}(F)$.
When explicitly stated, we use the same symbol for a representation or distribution of $G_{\mathrm{ad}}(F)$ and its inflation along $p$.

We identify the pinned dual groups with
\[
      \begin{aligned}
            \widehat{\rG_2}&=\rG_2(\C),& \widehat G&=\SO_8(\C),\\
            \widehat{G_{\mathrm{ad}}}&=\Spin_8(\C),& \widehat{G_{\mathrm{sim}}}&=\GSpin_8(\C).
      \end{aligned}
\]
The standard representation of $\Spin_8(\C)$ is denoted by $\Std_8$.

On the standard torus of $G_{\mathrm{sim}}$, write
\[
      t=\diag(t_1,t_2,t_3,t_4,\lambda/t_4,\lambda/t_3,\lambda/t_2,\lambda/t_1).
\]
The simple roots, which descend to $T_{G_{\mathrm{ad}}}$, are
\[
      \alpha_1(t)=t_1/t_2,\quad \alpha_2(t)=t_2/t_3,\quad
      \alpha_3(t)=t_3/t_4,\quad \alpha_4(t)=t_3t_4/\lambda.
\]
The central node is $\alpha_2$ and the outer nodes are $\alpha_1,\alpha_3,\alpha_4$; these $D_4$ roots are distinct from the $\rG_2$ roots $\alpha,\beta$ below.

The pinned diagram automorphisms of $G_{\mathrm{ad}}$ form $S_3$.
Fix triality $\theta$ of order three, so $A_3=\langle\theta\rangle$.
For an $S_3$-invariant irreducible $\pi'$ of $G_{\mathrm{ad}}(F)$, $\widetilde{\pi'}$ denotes its canonical extension to $G_{\mathrm{ad}}(F)\rtimes A_3$, as in Proposition~\ref{prop:canonical}.
Thus $\widetilde{\pi'}$ is a triality extension; the notation $\widetilde\pi$ in the similitude discussion denotes a representation of $G_{\mathrm{sim}}(F)$.
We distinguish the twisted coset $G_{\mathrm{ad}}(F) \rtimes \theta$ from the disconnected group $G_{\mathrm{ad}}(F)\rtimes\langle\theta\rangle$.
Let
\begin{align*}
      \widetilde{G_{\mathrm{ad}}} = \widetilde{\PGSO_8}
\end{align*}
denote the coset $G_{\mathrm{ad}}(F) \rtimes \theta$ with the twisted conjugation action of $G_{\mathrm{ad}}(F)$.

For the triality datum we put
\[
      \mathcal T=\operatorname{Tran}_{G_2}^{G_{\mathrm{ad}} \rtimes \theta, *},
      \qquad d(\theta)=|3|_F.
\]

For a parabolic $R$ of ${\rG_2}$, we suppress the superscript ${\rG_2}$ in $I_R^{\rG_2},r_R^{\rG_2},J_R^{\rG_2}$.
In particular, this defines $I_P,I_Q,J_P,J_Q$ used below.

Let $\alpha$ and $\beta$ be the short and long simple roots of ${\rG_2}$ relative to $B_{\rG_2}$.
For $\gamma\in\{\alpha,\beta\}$, let $P_\gamma$ be generated by $B_{\rG_2}$ and the negative root subgroup $U_{-\gamma}$.
Thus $\gamma$ is the simple root retained in its Levi.
Put $P=P_\alpha$ and $Q=P_\beta$, with Levi subgroups $M_P$ and $M_Q$.
Choose coordinates $t(t_1,t_2)$ on $T_{\rG_2}(F)$ by
\[
      \alpha(t(t_1,t_2))=t_1t_2^{-1},
      \qquad
      \beta(t(t_1,t_2))=t_2^2t_1^{-1}.
\]
Our identifications $M_P\simeq\GL_2$ and $M_Q\simeq\GL_2$ carry this element, respectively, to
\[
      \diag(t_1,t_2),\qquad
      \diag(t_2,t_1t_2^{-1}).
\]
All $\GL_2$ inducing representations on these Levis use these identifications.
In particular, the two determinant characters on $T_{\rG_2}(F)$ are $t_1t_2$ and $t_1$, respectively, and
\[
      \delta_P|_{M_P}=\nu_2^3,\qquad
      \delta_Q|_{M_Q}=\nu_2^5.
\]
Root lengths in $\widehat{\rG_2}$ refer to the dual root system: $\alpha^\vee$ is long and $\beta^\vee$ is short.
The names $P,Q$ always refer to the primal parabolics just fixed.

On either maximal Levi a twist $\tau\nu_2^s$ is interpreted using the displayed $\GL_2$ identification.
In particular, the notation $\pi(a,b)$ of Section~\ref{subsec:general-representations} always uses the intrinsic diagonal coordinates of $\GL_2$.

\subsection{Comparison with Qin's convention on $\rG_2$}\label{subsec:Qin-notation}

For comparison with \cite{Qin2026-AubertZelevinskyInvolutionG2HeckeAlgebras}, we record the conversion between ambient $\rG_2$ torus characters and the intrinsic $\GL_2$ notation of Section~\ref{subsec:general-representations}.

Let $I_T^\gamma(a\otimes b)$ denote normalized induction to $M_{P_\gamma}(F)$ of the character $t(t_1,t_2)\mapsto a(t_1)b(t_2)$ of $T_{\rG_2}(F)$.
The torus identifications in Section~\ref{subsec:pinned-groups} give
\[
      I_T^\alpha(a\otimes b)=\pi(a,b),
      \qquad I_T^\beta(a\otimes b)=\pi(ab,a).
\]
Indeed, on $M_Q$ the intrinsic character $(a,b)$ pulls back as $a(t_2)b(t_1/t_2)=b(t_1)(ab^{-1})(t_2)$.
Thus, intrinsic labels $(a,b)$ on $M_Q$ correspond to ambient labels $b\otimes ab^{-1}$.
The determinant characters computed above imply
\[
      \begin{aligned}
            I_T^\alpha(a\otimes b)\nu_2^s
               &=I_T^\alpha(\nu^sa\otimes\nu^sb),\\
            I_T^\beta(a\otimes b)\nu_2^s
               &=I_T^\beta(\nu^sa\otimes b).
      \end{aligned}
\]

In the induction notation $I_\gamma(s,\tau)$ and quotient notation $J_\gamma(s,\tau)$ used for these computations in \cite{Qin2026-AubertZelevinskyInvolutionG2HeckeAlgebras}, once $\tau$ has been expressed in intrinsic $\GL_2$ coordinates, the dictionary is
\[
      \begin{array}{c|cc}
            &\gamma=\alpha&\gamma=\beta\\ \hline
            I_\gamma(s,\tau)&I_P(\tau\nu_2^s)&I_Q(\tau\nu_2^s)\\
            J_\gamma(s,\tau)&J_P(\tau\nu_2^s)&J_Q(\tau\nu_2^s).
      \end{array}
\]
The second row applies to Langlands standard modules, and the Steinberg twist denoted $\delta(\chi)$ there is our $\St_2(\chi)$.

In particular, if $\chi^3=1$ and $\chi\ne1$, then ambient labels $\chi^{-1}\otimes\chi^{-1}$ on $M_Q$ give
\[
      I_T^\beta(\chi^{-1}\otimes\chi^{-1})
          =\pi(\chi^{-2},\chi^{-1})=\pi(\chi,\chi^{-1}).
\]
Throughout this paper the corresponding cubic quotient is therefore written $J_Q(\pi(\chi,\chi^{-1})\nu_2)$.
The symbol $\pi(a,b)$ always retains its intrinsic $\GL_2$ meaning.

\section{Representation theory of $\rG_2(F)$}

\subsection{Elliptic endoscopic groups of $\rG_2$}\label{subsec:elliptic-endoscopic-groups}

The elliptic endoscopic groups of $\rG_2$ are given by $\rG_2$, $\SO_4$ and $\PGL_3$.
Let $\sn$ be the spinor norm map
\begin{align}
      \sn \colon \SO_4(F) \to F^\times / F^{\times 2}.
\end{align}
We also have the determinant map
\begin{align}
      \det \colon \PGL_3(F) \to F^\times / F^{\times 3}.
\end{align}

Let $P_l$ (resp. $P_s$) be the parabolic subgroup of $\SO_4$ corresponding to the short root (resp. the long root) dual parabolic subgroup in $\widehat{\rG_2}(\C)$.

\subsection{\texorpdfstring{Arthur parameters for $\rG_2$ and their standard lifts}{Arthur parameters for G2 and their standard lifts}}\label{subsec:G2-parameters}

We apply Definition~\ref{def:parameters} to $\rG_2$.
Since $Z(\widehat{\rG_2})=1$, the component group used for ${\rG_2}$ is
\[
      \cS_\psi=\pi_0\Cent_{\widehat{\rG_2}}(\operatorname{Im}\psi).
\]
For a nontrivial Arthur $\SL_2$, its conjugacy class and centralizer are as follows:
\begin{center}
      \begin{tabular}{c|c}
            Arthur $\SL_2$ class & centralizer in $\widehat{\rG_2}$\\ \hline
            short-root & long-root $\SL_2$\\
            long-root & short-root $\SL_2$\\
            subregular & $S_3$\\
            principal & $1$
      \end{tabular}
\end{center}
The table is the list of nonzero nilpotent-orbit types and their $\SL_2$-centralizers in \cite{GanSavin2023-HoweDualityDichotomyExceptionalThetaCorrespondences}*{Section 2.7}.
The trivial Arthur factor gives the tempered case.
In the subregular case the connected Deligne factor maps trivially to the finite centralizer, so $\psi$ is co-tempered.
The possible subgroups $A=\psi(W_F)\subset S_3$, up to conjugacy, are $1$, $\mu_2$, $\mu_3$, and $S_3$.
For principal Arthur image the centralizer is trivial, so the Weil and Deligne factors are trivial.
For long-root and short-root parameters, the Deligne factor maps to the opposite-root $\SL_2$, either trivially or surjectively.
The former case is co-tempered; in the latter the Weil image lies in the center $\mu_2$ of this $\SL_2$.

\begin{definition}
      In the remaining cases, the image of the Deligne and Arthur $\SL_2$ appears in the standard representation as $S_2^D\boxtimes S_2^A$ or $S_2^A\boxtimes S_2^D$, respectively, where $S_n$ is the $n$-dimensional irreducible representation of $\SL_2(\C)$.
      We call these parameters mixed parameters.
\end{definition}

Fix the standard embedding $\iota = \iota^{\PGSO_8}_{\rG_2} : \widehat{\rG_2}\hookrightarrow\widehat{G_{\mathrm{ad}}}=\Spin_8(\C)$.
Put $\psi'=\iota\circ\psi$.
The dual maps of $q$ and $p$ are $\widehat q:\widehat{G_{\mathrm{ad}}}\to\widehat G$ and $\widehat p:\widehat{G_{\mathrm{ad}}}\hookrightarrow\widehat{G_{\mathrm{sim}}}$.
We write $\psi_0=\widehat q\circ\psi'=\Std_8\circ\psi'$ for the orthogonal parameter and $\widetilde\psi=\widehat p\circ\psi'$ for its similitude lift.
For the two commuting root $\SL_2$ factors, the standard representation restricts as
\[
      \Std_8=(S_{2,s}\boxtimes S_{2,l})\oplus S_{3,s}\oplus\mathbf1.
\]
Here $s$ and $l$ specify the short and long factors in $\widehat{\rG_2}$.

\begin{corollary}\label{cor:mixed-parameters}
      For a quadratic character $\chi$, including $\chi=1$, the long-root Arthur parameter with Deligne parameter $\chi\boxtimes S_2$ in the short-root centralizer has standard lift
      \begin{align*}
            \Std_8\circ\psi'_D&=\chi\boxtimes S_2^D\boxtimes S_2^A
                                  \oplus S_3^D\oplus\mathbf1,\\
            \Std_8\circ\phi_{\psi'_D}&=
             \chi\nu^{1/2}\boxtimes S_2^D\oplus\chi\nu^{-1/2}\boxtimes S_2^D
                                  \oplus S_3^D\oplus\mathbf1.
      \end{align*}
      For the swapped, short-root parameter the corresponding expressions are
      \begin{align*}
            \Std_8\circ\psi'_A&=\chi\boxtimes S_2^D\boxtimes S_2^A
                                  \oplus S_3^A\oplus\mathbf1,\\
            \Std_8\circ\phi_{\psi'_A}&=
             \chi\nu^{1/2}\boxtimes S_2^D\oplus\chi\nu^{-1/2}\boxtimes S_2^D
                                  \oplus\nu\oplus\mathbf1\oplus\nu^{-1}\oplus\mathbf1.
      \end{align*}
\end{corollary}

\begin{proof}
      Use the restriction of $\Std_8$ above and $\chi^2=1$, and then apply the block formula~\eqref{eq:parameter-block}.
      Swapping the factors interchanges the $S_3^D$ and $S_3^A$ summands.
\end{proof}

Subregular parameters factor through
\begin{align*}
      L_F \times \SL_2^A \to W_F \times \SL_2^A \xrightarrow{\sigma \times S_3^A} S_3 \times \SO_3 \hookrightarrow \rG_2.
\end{align*}
The action of $S_3 \times \SO_3$ on the standard representation of $\rG_2$ is given by
\begin{align*}
      \rho \boxtimes S_3^A \oplus \epsilon \boxtimes S_1^A,
\end{align*}
where $\rho$ is the standard 2-dimensional representation of $S_3$ and $\epsilon$ is the sign character of $S_3$.

\begin{lemma}\label{lemma:subregular-standard}
      The lift $\psi'_{\sigma} = \iota \circ \psi_{\sigma}$ of a subregular A-parameter
      \begin{align*}
            \psi_{\sigma} \colon L_F \times \SL_2^A \to W_F \times \SL_2^A \xrightarrow{\sigma \times S_3^A} S_3 \times \SO_3 \hookrightarrow \rG_2.
      \end{align*}
      to $\PGSO_8$ has the standard parameter
      \begin{align*}
            \Std_8 \circ \psi'_{\sigma} = \rho \circ \sigma \boxtimes S_3^A \oplus \det \rho \circ \sigma \oplus \mathbf{1}.
      \end{align*}
      The swapped parameter is a tempered $L$-parameter
      \begin{align*}
            \rho \circ \sigma \boxtimes S_3^D \oplus \det \rho \circ \sigma \oplus \mathbf{1}.
      \end{align*}
\end{lemma}

\subsection{A rank-one Aubert calculation}

\begin{lemma}\label{lem:rank-one}
      Let $R$ be $P$ or $Q$, and let $\tau$ be unitary supercuspidal on its $\GL_2$ Levi.
      Suppose the standard module $I_R(\tau\nu_2^s)$, $s>0$, has length two, with a square-integrable subrepresentation $\delta$ and Langlands quotient $J$.
      Then
      \[
            [J]=[I_R(\tau\nu_2^s)]-[\delta],\qquad
            \mathcal D_{\rG_2}[\delta]=-[J],\qquad\widehat\delta=J.
      \]
\end{lemma}

\begin{proof}
      The geometric lemma gives the two cuspidal exponents $\tau\nu_2^s$ and its Weyl conjugate in $r_R I_R(\tau\nu_2^s)$, each with multiplicity one.
      Frobenius reciprocity supplies the first exponent in $r_R\delta$; the square-integrability criterion excludes the opposite one.
      Thus $r_R\delta=\tau\nu_2^s$.
      Its Jacquet modules for the torus and for the nonconjugate maximal Levi vanish, since its supercuspidal support is in $M_R$.
      Formula~\eqref{eq:Aubert-operator} now gives $\mathcal D_{\rG_2}[\delta]=[\delta]-[I_R(\tau\nu_2^s)]=-[J]$.
      The support Levi has semisimple rank one, so $\epsilon_{\rG_2}(\delta)=-1$ and Proposition~\ref{prop:Aubert} gives the last assertion.
      The exact sequences used below are supplied by \cite{GanSavin2023-HoweDualityDichotomyExceptionalThetaCorrespondences}*{ Propositions 3.1(i) and 3.2(i)}.
\end{proof}

\section{Representation theory of $\widetilde{\PGSO_8}(F)$}

\subsection{Similitude packets and the distinguished tempered lift}\label{subsec:similitude-preliminaries}

Let $\Sigma_0=\langle\theta_0\rangle$, where $\theta_0$ is the standard pinned outer involution of $G$, extended to $G_{\mathrm{sim}}$; we use the corresponding involutions on the dual groups.
Bars on $\Irr$ and packet symbols denote passage to $\Sigma_0$-orbits, and brackets denote these orbits on representations or parameter conjugacy classes, using the dual involution for parameters.

We use the infinitesimal character of Section~\ref{subsec:parameters-and-infinitesimal-characters}, in particular $\lambda_{\widetilde\psi}(w)=\widetilde\psi(w,d_w,d_w)$.
Fix the local Langlands correspondence modulo $\Sigma_0$ with the choices of \cite{Xu2024-ArthurPacketsQuasisplitGSp2nGO2n}*{Section 5}.
For $[\widetilde\pi]\mapsto[\widetilde\varphi_{\widetilde\pi}]$, put $[\lambda_{\widetilde\pi}] =[\lambda_{\widetilde\varphi_{\widetilde\pi}}]$ and define
\begin{equation}\label{eq:Irr-infinitesimal-character}
      \overline{\Irr}(G_{\mathrm{sim}}(F))_{\widetilde\lambda}
      =\left\{[\widetilde\pi]\in\Irr(G_{\mathrm{sim}}(F))/\Sigma_0:
            [\lambda_{\widetilde\pi}]=[\widetilde\lambda]\right\}.
\end{equation}
Thus the subscript fixes the infinitesimal character in the sense of \cite{Xu2024-ArthurPacketsQuasisplitGSp2nGO2n}*{Section 3, pp.~10353--10356}.

For $\widetilde\xi=\widetilde\lambda$ or $\widetilde\psi$, write
\[
      X^{\Sigma_0}(\widetilde\xi)=
      \left\{\omega\in\Hom(G_{\mathrm{sim}}(F)/G(F),\C^\times):
            [\widetilde\xi\otimes\omega]=[\widetilde\xi]\right\}.
\]
Here $\omega$ is a smooth character, identified with a character of $W_F$ by local class field theory, and the twist uses the central $\C^\times$ in $\widehat{G_{\mathrm{sim}}}\to\widehat G$.

\begin{proposition}\label{prop:Xu-packets}
      Fix a central character of $G_{\mathrm{sim}}(F)$ extending that of a tempered packet of $G(F)$.
      The refined similitude packets of \cite{Xu2018-LPacketsQuasisplitGSp2nGO2n}*{Theorem 4.6} have stable sums of characters and satisfy the character relations of that theorem.
      Their choice is unique only up to the character group
      \[
            X=\Hom\bigl(G_{\mathrm{sim}}(F)/(Z_{G_{\mathrm{sim}}}(F)G(F)),\C^\times\bigr).
      \]
      For Arthur parameters, the corresponding existence, stability, and character identities are supplied by \cite{Xu2024-ArthurPacketsQuasisplitGSp2nGO2n}*{Theorem 6.4}, with its central-character, outer-conjugacy, and normalization conventions.
\end{proposition}

\begin{proof}
      The tempered assertion follows from \cite{Xu2018-LPacketsQuasisplitGSp2nGO2n}*{Theorem 4.6(1)--(3)}, and the twisting ambiguity is described in \cite{Xu2018-LPacketsQuasisplitGSp2nGO2n}*{Remark 4.7}.
      The Arthur assertion is \cite{Xu2024-ArthurPacketsQuasisplitGSp2nGO2n}*{Theorem 6.4(1)--(3)}.
\end{proof}

For $G_{\mathrm{ad}}$, we impose the trivial scalar central character on $G_{\mathrm{sim}}(F)$.
If $\pi'$ is a representation of $G_{\mathrm{ad}}(F)$, the similitude parameter of $p^*\pi'$ therefore factors through $\Spin_8(\C)\subset\GSpin_8(\C)$.
Write $\phi_{\pi'}$ for a representative of this parameter modulo the standard outer involution.
The linear parameter $\Std_8\circ\phi_{\pi'}$ is independent of this choice up to isomorphism.
In the following criterion, $\gamma(s,r,\psi_F)$ is the Artin gamma factor of the linear Weil--Deligne representation $r$; the parameter of $s_1\pi'$ is obtained separately in the same way.
We use the following distinguished tempered packet.

\begin{theorem}\label{thm:distinguished}
      Let $s_1$ be the pinned involution interchanging the standard and one half-spin representation of $\Spin_8$.
      A tempered irreducible $\pi'$ of $G_{\mathrm{ad}}(F)$ satisfying
      \[
            \gamma(0,\Std_8\circ\phi_{\pi'},\psi_F)=
            \gamma(0,\Std_8\circ\phi_{s_1\pi'},\psi_F)=0
      \]
      is the composite theta lift of a unique tempered representation of ${\rG_2}(F)$.
      The map taking $\pi'$ to the parameter of this representation is surjective onto the tempered parameters of ${\rG_2}$.
      Its fibers are distinguished Xu packets, denoted $\Pi^X_{\iota\circ\varphi}$.
      Every member is invariant under all pinned diagram automorphisms $S_3$ of $G_{\mathrm{ad}}$.
\end{theorem}

\begin{proof}
      The characterization and the assertion about fibers are \cite{TakanashiThesisLocal}*{Proposition 3.58 and Definition 3.59, pp.~25--26}.
      The $S_3$-invariance follows from \cite{TakanashiThesisLocal}*{Proposition 6.6}.
      The underlying tempered correspondence for ${\rG_2}$ is \cite{GanSavin2023-LocalLanglandsConjectureG2}; the similitude packets are those of Proposition~\ref{prop:Xu-packets}.
\end{proof}

\subsection{Restriction and compatibility with Arthur packets}

We use the maps $j,p,q$ of Section~\ref{subsec:pinned-groups} and the positive Aubert functor $D_A$ of Section~\ref{subsec:aubert-involution-and-endoscopic-transfer}.
For $G_{\mathrm{ad}}(F)\rtimes S_3$, the same symbol denotes the positive disconnected Aubert functor of \cite{ChengKaletha2026-AubertBernsteinInvolutionsDisconnectedGroups}*{Definition 3.1.9}.

\begin{lemma}\label{lem:Aubert-restriction}
      Restriction from $G_{\mathrm{sim}}(F)$ to $G(F)$, pullback along $q$, and restriction from $G_{\mathrm{ad}}(F)\rtimes S_3$ to $G_{\mathrm{ad}}(F)$ commute with $D_A$.
\end{lemma}

\begin{proof}
      Corresponding standard parabolics $\widetilde R$ and $R$ of $G_{\mathrm{sim}}$ and $G$ have the same unipotent radical and modulus character, and $G_{\mathrm{sim}}(F)=\widetilde R(F)G(F)$.
      Thus restriction commutes with normalized induction, Jacquet modules, and the maps in the Aubert complex.
      The relative semisimple ranks agree, so taking its cohomology gives the assertion for $D_A$.
      Inflation from $G_{\mathrm{ad}}$ to $G_{\mathrm{sim}}$ also identifies these complexes, proving the assertion for $q$.
      The last assertion is \cite{ChengKaletha2026-AubertBernsteinInvolutionsDisconnectedGroups}*{ Lemma 3.1.10}.
\end{proof}

\begin{proposition}\label{prop:standard-Aubert}
      Let $\psi_0=\Std_8\circ\psi'$.
      If $\operatorname{JH}(q^*\pi')\subset\Pi_{\psi_0}(G)$, then
      \[
            \operatorname{JH}(q^*\widehat{\pi'})
                  \subset\Pi_{\widehat{\psi_0}}(G),\qquad
            \widehat{\psi_0}=\Std_8\circ\widehat{\psi'}.
      \]
      Packets of $G$ are taken modulo its standard outer involution.
\end{proposition}

\begin{proof}
      For classical Arthur packets, $\widehat{\Pi_{\psi_0}}=\Pi_{\widehat{\psi_0}}$ by \cite{AtobeGanEtAl2025-LocalIntertwiningRelationsCoTemperedPackets}*{ Lemma 4.4.4(1)}, whose character-relation hypotheses follow from Arthur's local classification, as recalled in \cite{Xu2017-MoeglinParametrizationArthurPackets}*{Theorems 4.1--4.2}.
      In the even orthogonal case the cited lemma is formulated for the full orthogonal group $\mathrm O_8(F)$.
      Restricting its packet members to $G(F)=\SO_8(F)$ recovers the classical packets modulo $\Sigma_0$; see \cite{AtobeGanEtAl2025-LocalIntertwiningRelationsCoTemperedPackets}*{ Remark 1.6.2(1)}.
      Positive Aubert duality commutes with restriction to the identity component by \cite{ChengKaletha2026-AubertBernsteinInvolutionsDisconnectedGroups}*{ Lemma 3.1.10}, so the same packet identity holds in these outer-involution orbits.
      Apply Lemma~\ref{lem:Aubert-restriction} to $q^*\pi'$.
      The parameter equality follows by interchanging the two $\SL_2$ factors.
\end{proof}

\begin{proposition}\label{prop:Xu-Aubert}
      Let $\widetilde\psi$ be the $G_{\mathrm{sim}}$ parameter induced by $\psi'$, and put $\widetilde\lambda=\lambda_{\widetilde\psi}$.
      In Xu's notation, suppose
      \begin{equation}\label{eq:Xu-criterion}
            X^{\Sigma_0}(\widetilde\lambda)
                  =X^{\Sigma_0}(\widetilde\psi).
      \end{equation}
      Then the positive Aubert duals of the Xu packet for $\psi'$ form the Xu packet for $\widehat{\psi'}$, with the same infinitesimal and central characters.
      Condition~\eqref{eq:Xu-criterion} holds if $\psi'$ is tempered or co-tempered.

      Furthermore, in this case the restriction map induces a bijection
      \begin{align*}
            \overline{\Pi}_{\widetilde{\psi}}(G_{\mathrm{sim}}) &\longrightarrow \overline{\Pi}_{\psi_0}(G)/G_{\mathrm{sim}}(F).
      \end{align*}
\end{proposition}

\begin{proof}
      By \cite{Xu2024-ArthurPacketsQuasisplitGSp2nGO2n}*{Theorem 6.5 and the preceding construction}, the packet is
      \[
            \overline\Pi_{\widetilde\psi}
            =\left\{[\widetilde\pi]\in
                  \overline{\operatorname{Irr}}(G_{\mathrm{sim}}(F))_{\widetilde\lambda}:
                  \operatorname{JH}(\widetilde\pi|_{G(F)})
                                    \subset\overline\Pi_{\psi_0}\right\}.
      \]
      The notation is as in Section~\ref{subsec:similitude-preliminaries}.
      Aubert duality preserves cuspidal support, hence the infinitesimal character by \cite{Xu2024-ArthurPacketsQuasisplitGSp2nGO2n}*{Proposition 5.1 and Corollary 5.5}; it also preserves the central character.
      Moreover,
      \[
            \lambda_{\widehat{\widetilde\psi}}=\widetilde\lambda,
            \qquad
            X^{\Sigma_0}(\widehat{\widetilde\psi})
                  =X^{\Sigma_0}(\widetilde\psi),
      \]
      since swapping the two $\SL_2$ factors changes neither the diagonal parameter nor the image.
      Proposition~\ref{prop:standard-Aubert} gives one inclusion between the packets; involutivity gives the reverse.
      For a tempered or co-tempered parameter the diagonal parameter $\widetilde\psi_d(w,x)=\widetilde\psi(w,x,x)$ has the same image as $\widetilde\psi$ and is the open-orbit parameter for $\widetilde\lambda$.
      By \cite{Xu2024-ArthurPacketsQuasisplitGSp2nGO2n}*{Section 3, p.~10356}, we have $X^{\Sigma_0}(\widetilde\psi_d)=X^{\Sigma_0}(\widetilde\lambda)$, proving \eqref{eq:Xu-criterion}.
\end{proof}

\section{Tempered twisted endoscopy for $\widetilde{\PGSO_8}$ and Aubert duality}\label{sec:twisted-endoscopy}

\subsection{Tempered character identities and stability on proper Levi subgroups}

For a tempered parameter $\varphi$ of $\rG_2$, write $\Pi_\varphi$ for its Gan--Savin packet and $\rho_{\varphi,\pi}\in\Irr(\cS_\varphi)$ for the enhancement in the fixed Whittaker normalization constructed in \cite{TakanashiEndoscopyLocal}*{Theorem 31.1}.
Set
\[
      m_\pi=\dim\rho_{\varphi,\pi},\qquad
      S\Theta_\varphi=\sum_{\pi\in\Pi_\varphi}m_\pi\Theta_\pi.
\]
This is the stable tempered distribution used below.

\begin{theorem}\label{thm:tempered-input}
      For each tempered parameter $\varphi$ considered below, one has
      \begin{equation}\label{eq:tempered-triality}
            \mathcal T(S\Theta_\varphi)=
            \sum_{\pi'\in\Pi^X_{\iota\circ\varphi}}\Theta_{\widetilde{\pi'}}.
      \end{equation}
      More generally, we also have the twisted endoscopic character identities
      \begin{align*}
            \sum_{\pi'\in\Pi^X_{\iota\circ\varphi}} \langle s, \pi' \rangle \Theta_{\widetilde{\pi'}}(f^{\widetilde{\PGSO_8}})
            = S\Theta^{H_s}_{\varphi^s}(\Trans_{H_s}^{\widetilde{\PGSO_8}} f^{\widetilde{\PGSO_8}}),
      \end{align*}
      for all $s\in\cS_{\iota_{\rG_2}^{\PGSO_8}\circ\varphi}$, where $\varphi^s$ is the parameter of the endoscopic group $H_s$ corresponding to $s$.

      The tempered endoscopic identities for $\rG_2$ used in Sections~\ref{sec:subregular}--\ref{sec:mixed-short} hold with the same normalization.
\end{theorem}

\begin{proof}
      This is the non-archimedean case of \cite{TakanashiEndoscopyLocal}*{Theorem 31.1}, applied to the distinguished packets of Theorem~\ref{thm:distinguished}.
\end{proof}

The proper Levis on the $D_4$ side have the following additional stability property.

\begin{lemma}\label{lem:proper-Levi-stable}
      Every finite-length character of a proper Levi subgroup of $G_{\mathrm{sim}}$ or $G_{\mathrm{ad}}$ is stable.
      Normalized induction from a proper Levi therefore preserves stability without any condition on the inducing finite-length representation.
\end{lemma}

\begin{proof}
      Such a Levi is a quotient $L/Z$ of a product of general linear groups by a central split torus.
      Indeed, up to the outer involution, proper Levis of $G_{\mathrm{sim}}$ are products of general linear groups and either a split torus or $\GSO_{2m}$ with $m\leq3$, and
      \[
            \GSO_4\simeq(\GL_2\times\GL_2)/\{(zI_2,z^{-1}I_2)\},\qquad
            \GSO_6\simeq(\GL_4\times\GL_1)/\{(zI_4,z^{-2})\}.
      \]
      Passing to $G_{\mathrm{ad}}$ quotients by another central split torus.
      Lemma~\ref{lem:central-quotient-stability} therefore gives stability of every finite-length character of such a Levi.
      Apply Lemma~\ref{lem:induction-transfer} for the last assertion.
\end{proof}

We fix the inherited Levi embeddings for the later calculations.
For $R=P,Q$ with $\GL_2$ Levi $M_R$, let $M_{R'}$ be the $\theta$-stable Levi of $G_{\mathrm{ad}}$ whose dual group is
\[
      \widehat M_{R'}=\Cent(\iota(Z(\widehat M_R)^\circ),\Spin_8(\C)).
\]
Take the standard parabolic $R'=M_{R'}N_{R'}$ with this Levi.
Restriction of $\iota$ gives the inherited Levi embedding; denote its transfer by $\mathcal T_{M_R}^{M_{R'},\theta}$.
The derived types of $M_{P'},M_{Q'}$ are $A_1,A_1^3$, respectively.
Thus Lemma~\ref{lem:induction-transfer} gives, for stable virtual characters $S$,
\[
      \mathcal T(I_R S)=I_{R'}^{G_{\mathrm{ad}}}
      \mathcal T_{M_R}^{M_{R'} \rtimes \theta}(S).
\]

The parabolic subgroup $Q' = P_{\alpha_1, \alpha_3, \alpha_4}$ of $\PGSO_8$ has Levi subgroup $M_{Q'} = M_{\alpha_1, \alpha_3, \alpha_4} = (\GL_2^3 \times \GL_1)/ \GL_1^3$, where the embedding is given by $(z_1, z_2, z_3) \mapsto (z_1, z_2, z_3, z_1 z_2 z_3)$.
The action of the triality group $S_3$ on $M_{Q'}$ is given by permuting the three $\GL_2$-factors.
For any representation $\sigma$ of $M_Q \cong \GL_2$, we define a representation of $M_{Q'}$ by
\begin{align*}
      \sigma \boxtimes \sigma \boxtimes \sigma \boxtimes \omega_{\sigma}^{-1}
\end{align*}
and denote it by $(\sigma)_{M_{Q'}}$.
Then, the canonical extension to $M_{Q'} \rtimes \theta$ is equal to $\Trans_{M_{Q}}^{M_{Q'} \rtimes \theta}(\sigma)$.
We have $M_{P'} \cong \GL_2 \times \GL_1^3/\GL_1$ where the embedding is given by $z \mapsto (zI_2, z, z, z)$.
We use similar notation $(\sigma)_{M_{P'}}$ for the representation of $M_{P'}$ defined by
\begin{align*}
      \sigma \omega_{\sigma} \boxtimes \omega_{\sigma}^{-1} \boxtimes \omega_{\sigma}^{-1} \boxtimes \omega_{\sigma}^{-1}.
\end{align*}
For the standard Borel subgroup $B_{\PGSO_8}$ of $\PGSO_8$, a triality-stable representation of the maximal torus $T_{\PGSO_8}$ corresponds to a unique representation of $T_{\rG_2}$ via
\begin{align*}
      (T_{\PGSO_8})_{\theta} \cong T_{\rG_2}.
\end{align*}
For a character $\chi$ of $T_{\rG_2}$, we denote the corresponding character of $T_{\PGSO_8}$ by $(\chi)_{T_{\PGSO_8}}$.

From the Langlands classification, we obtain the following result.
\begin{lemma}\label{lem:triality-invariant-non-tempered}
      Let $(R', \sigma')$ be a pair consisting of a proper standard parabolic subgroup $R' = M_{R'} N_{R'}$ of $\PGSO_8$ and an essentially tempered representation $\sigma'$ of $M_{R'}$ such that the central character of $\sigma'$ is positive with respect to $R'$.
      Then, the corresponding Langlands quotient $J_{R'}(\sigma')$ of $\PGSO_8$ is triality-invariant if and only if $R' \in \{B_{\PGSO_8}, P', Q' \}$ and $\sigma'$ is triality-invariant.

      In this case, there exists a unique pair of a standard parabolic subgroup $R$ of $\rG_2$ and an essentially tempered representation $\sigma$ of $M_R$ such that $R'$ corresponds to $R$ and $\sigma' = (\sigma)_{M_{R'}}$ corresponds to $\sigma$.
\end{lemma}

\subsection{Triality signs and Aubert transfer}

Specializing \eqref{eq:twisted-Aubert-operator} gives
\[
      \mathcal D_{G_{\mathrm{ad}}}^\theta=
       \sum_{\substack{R=MN\supset B_{G_{\mathrm{ad}}}\\\theta(R)=R}}
         (-1)^{\dim A_M^\theta}I_R^{G_{\mathrm{ad}}}r_R^{G_{\mathrm{ad}}}.
\]
The following proposition identifies its support sign on canonical triality extensions.

\begin{proposition}\label{prop:twisted-sign}
      Let $\pi'$ be $S_3$-invariant and irreducible.
      Then
      \begin{equation}\label{eq:twisted-sign}
            \mathcal D_{G_{\mathrm{ad}}}^\theta\Theta_{\widetilde{\pi'}}
                  =\epsilon_{G_{\mathrm{ad}}}(\pi')\Theta_{\widetilde{\widehat{\pi'}}}
                  \quad\text{on }G_{\mathrm{ad}}(F)\theta.
      \end{equation}
\end{proposition}

\begin{proof}
      The Aubert dual commutes with restriction to $G_{\mathrm{ad}}(F)$ and preserves irreducibility, up to sign, by \cite{ChengKaletha2026-AubertBernsteinInvolutionsDisconnectedGroups}*{ Lemmas 3.1.10--3.1.11}.
      Apply it to an $S_3$-extension of $\pi'$.
      The resulting representation is an $S_3$-extension of $\widehat{\pi'}$, so Proposition~\ref{prop:canonical} identifies its restriction to $A_3$ with the canonical extension.
      The character formula and sign comparison in \cite{ChengKaletha2026-AubertBernsteinInvolutionsDisconnectedGroups}*{ Proposition 3.2.1 and Corollary 3.2.4} give the scalar $\epsilon_{G_{\mathrm{ad}}}(\pi')$.
      Note that the orientation twist in \cite{ChengKaletha2026-AubertBernsteinInvolutionsDisconnectedGroups}*{Corollary 3.3.2} is trivial on $\theta$.
\end{proof}

\begin{theorem}\label{thm:signed-transfer}
      Suppose a stable virtual character $S=\sum_i a_i\Theta_{\pi_i}$ of ${\rG_2}$ satisfies $\mathcal T(S)=\sum_j b_j\Theta_{\widetilde{\pi'_j}}$, with the $\pi'_j$ irreducible and $S_3$-invariant.

      Then, we have
      \begin{equation}\label{eq:signed-transfer}
            \mathcal T\left(\sum_i a_i\epsilon_{\rG_2}(\pi_i)
                                    \Theta_{\widehat\pi_i}\right)
             =\sum_j b_j \epsilon_{G_{\mathrm{ad}}}(\pi'_j)
                  \Theta_{\widetilde{\widehat{\pi'_j}}}.
      \end{equation}
      The distribution inside $\mathcal T$ is stable.
      The ordinary analogue for $\rG_2$ follows from the compatibility formula in Section~\ref{subsec:aubert-involution-and-endoscopic-transfer}.
\end{theorem}

\begin{proof}
      Apply $\mathcal D_{G_{\mathrm{ad}}}^\theta$ to the assumed equality, move it across $\mathcal T$ by Proposition~\ref{prop:Xu-transfer}, and use Propositions~\ref{prop:Aubert} and \ref{prop:twisted-sign} term by term.
      Stability is the left vertical assertion of \cite{Xu2017-MoeglinParametrizationArthurPackets}*{diagram (A.1)}.
      For ordinary transfer the same argument uses Hiraga's formula.
\end{proof}

\begin{corollary}\label{cor:cotempered}

      For a co-tempered parameter $\psi$, put $\varphi=\widehat\psi$ and define
      \[
            S\Theta_\psi^{\mathrm{Aub}}=\mathcal D_{\rG_2} S\Theta_\varphi,
            \qquad T_\psi^{\mathrm{Aub}}=
                  \mathcal D_{G_{\mathrm{ad}}}^\theta\sum_{\pi'\in\Pi^X_{\iota\circ\varphi}}
                                                 \Theta_{\widetilde{\pi'}}.
      \]
      Then $S\Theta_\psi^{\mathrm{Aub}}$ is stable and $\mathcal T(S\Theta_\psi^{\mathrm{Aub}})=T_\psi^{\mathrm{Aub}}$.
      Its irreducible terms are the Aubert duals of those of $S\Theta_\varphi$, with the coefficients prescribed in \eqref{eq:signed-transfer}.
      On $G_{\mathrm{ad}}$ these duals form the distinguished Xu packet
      \[
            \Pi^X_{\iota\circ\psi}
                  =\{\widehat{\pi'}:\pi'\in\Pi^X_{\iota\circ\varphi}\}.
      \]
\end{corollary}

\begin{proof}
      Apply Theorem~\ref{thm:signed-transfer} to \eqref{eq:tempered-triality} and Proposition~\ref{prop:Xu-Aubert}.
      Lemma~\ref{lem:Aubert-restriction} preserves the $S_3$-extensions, whose restrictions to $A_3$ are canonical by Proposition~\ref{prop:canonical}.
\end{proof}

The support signs are
\[
      \begin{array}{c|c|c}
            \text{group}&\text{semisimple type of the supercuspidal support Levi}&\text{support sign}\\\hline
            {\rG_2}&\text{torus}&+1\\
            {\rG_2}&A_1&-1\\
            {\rG_2}&\rG_2&+1\\
            G_{\mathrm{ad}}&\text{torus, }A_1^2\text{, or }A_2&+1\\
            G_{\mathrm{ad}}&A_1\text{, }A_3\text{, or }A_1^3&-1\\
            G_{\mathrm{ad}}&D_4&+1
      \end{array}
\]

\section{The principal parameter}\label{sec:principal-parameter}

\begin{proposition}\label{prop:principal-ecr}
      For the principal Arthur parameter, the Aubert-dual set is $\{\mathbf1_{\rG_2}\}$, where $\mathbf1_{\rG_2}=J_{B_{\rG_2}}^{\rG_2}(\delta_{B_{\rG_2}}^{1/2})$, and
      \[
            \mathcal T(\Theta_{\mathbf1_{\rG_2}})=\Theta_{\mathbf1_{\PGSO_8}},
      \]
      where the trivial representation on the right has the trivial triality action.
\end{proposition}

\begin{proof}
      The swapped parameter has packet $\{\St_{\rG_2}\}$.
      Its lift is principal in $\Spin_8$: its standard representation is $S_7\oplus\mathbf1$.
      The distinguished tempered lift is $\St_{G_{\mathrm{ad}}}$, and Theorem~\ref{thm:tempered-input} gives $\mathcal T(\Theta_{\St_{\rG_2}})=\Theta_{\widetilde{\St_{G_{\mathrm{ad}}}}}$.
      The positive Aubert dual of the Steinberg representation is the trivial representation.
      Both Steinberg representations have toral supercuspidal support, and the ranks are two and four, so both signs are $+1$.
      Theorem~\ref{thm:signed-transfer} proves the assertion.
\end{proof}

\section{Subregular parameters}\label{sec:subregular}

\subsection{The tempered subregular packet}\label{sec:tempered-subregular}

Let $\psi$ be a subregular Arthur parameter of ${\rG_2}$.
Then, the swapped parameter $\phi = \widehat\psi$ is tempered.
We use the description in Section \ref{subsec:G2-parameters} and write $\phi = \phi_\sigma$ for a homomorphism
\begin{align*}
      \sigma \colon W_F \to S_3.
\end{align*}
Its image is either $1$, $\mu_2$, $\mu_3$, or $S_3$.
In the second and third cases, let $\chi$ be the corresponding nontrivial quadratic or cubic character of $W_F$.

\begin{theorem}\label{thm:tempered-subregular-elements}
      We have the following tempered subregular $L$-packets of ${\rG_2}$.
      \begin{enumerate}
            \item If $\sigma(W_F) = \{1\}$, then the $L$-packet $\Pi_{\phi}$ consists of three representations $\pi_{\gen}[1], \pi_{\deg}[1], \pi_{\mathrm{sc}}[1]$ corresponding to the three irreducible representations $\mathbf{1}, \rho, \epsilon$ of $S_3$.
            The representation $\pi_{\gen}[1]$ is generic, $\pi_{\deg}[1]$ is non-generic, and $\pi_{\mathrm{sc}}[1]$ is supercuspidal.
            \item If $\sigma(W_F) = \mu_2$, then the $L$-packet $\Pi_{\phi}$ consists of two representations $\pi_{\gen}[\chi], \pi_{\mathrm{sc}}[\chi]$ corresponding to the two irreducible representations $\mathbf{1}, \epsilon$ of $\mu_2$.
            The representation $\pi_{\gen}[\chi]$ is generic, and $\pi_{\mathrm{sc}}[\chi]$ is non-generic.
            \item If $\sigma(W_F) = \mu_3$, then the $L$-packet $\Pi_{\phi}$ consists of three representations $\pi_{\gen}[\chi], \pi_{\mathrm{sc}, 1}[\chi], \pi_{\mathrm{sc}, 2}[\chi]$.
            The representation $\pi_{\gen}[\chi] = \pi_{\gen}[\chi^2]$ is generic, and $\pi_{\mathrm{sc}, 1}[\chi], \pi_{\mathrm{sc}, 2}[\chi]$ are supercuspidal.
            \item If $\sigma(W_F) = S_3$, then the $L$-packet $\Pi_{\phi}$ consists of a single generic representation $\pi_{\gen}[\rho \circ \sigma]$.
      \end{enumerate}

      In Cases $(1), (2), (3)$ all the non-supercuspidal representations in the $L$-packets have the cuspidal support in the maximal torus of ${\rG_2}$.
      In Case $(4)$, the generic representation $\pi_{\gen}[\rho \circ \sigma]$ has the cuspidal support in the proper Levi subgroup $M_Q$ corresponding to the long root of ${\rG_2}$.
\end{theorem}

\begin{proposition}\label{prop:tempered-subregular}
      The stable distribution $S\Theta_{\phi_{\sigma}}$ of the tempered subregular $L$-packet $\Pi_{\phi_{\sigma}}$ is given by
      \begin{align*}
            S\Theta_{\phi_{\sigma}} = \sum_{\pi \in \Pi_{\phi}} m_{\pi} \Theta_{\pi},
      \end{align*}
      where $m_{\pi}$ is the dimension of the corresponding irreducible representation of $\cS_{\phi_{\sigma}} = \Cent_{S_3}(\sigma(W_F))$.

      We have the standard endoscopic character identities:
      \begin{align*}
            \Trans^{\rG_2, *}_{\SO_4}(\Theta_{\St_{\SO_4}}) &= \Theta_{\pi_{\gen}[1]} - \Theta_{\pi_{\mathrm{sc}}[1]}, \\
            \Trans^{\rG_2, *}_{\SO_4}(\Theta_{\St_{\SO_4} \otimes \chi \circ \sn}) &= \Theta_{\pi_{\gen}[\chi]} - \Theta_{\pi_{\mathrm{sc}}[\chi]}, &&(\chi^2=1,\ \chi\ne1) \\
            \Trans^{\rG_2, *}_{\PGL_3}(\Theta_{\St_{\PGL_3}}) &= \Theta_{\pi_{\gen}[1]} - \Theta_{\pi_{\deg}[1]} + \Theta_{\pi_{\mathrm{sc}}[1]}.
      \end{align*}

      For a non-trivial cubic character $\chi$, we can choose a primitive root of unity $\zeta$ so that we have
      \begin{align*}
            \Trans^{\rG_2, *}_{\PGL_3}(\Theta_{\St_{\PGL_3} \otimes \chi \circ \det}) &=
            \Theta_{\pi_{\gen}[\chi]} +\zeta \Theta_{\pi_{\mathrm{sc}, 1}[\chi]} + \zeta^2 \Theta_{\pi_{\mathrm{sc}, 2}[\chi]}, \\
            \Trans^{\rG_2, *}_{\PGL_3}(\Theta_{\St_{\PGL_3} \otimes \chi^2 \circ \det}) &= \Theta_{\pi_{\gen}[\chi]} + \zeta^2 \Theta_{\pi_{\mathrm{sc}, 1}[\chi]} + \zeta \Theta_{\pi_{\mathrm{sc}, 2}[\chi]}.
      \end{align*}
\end{proposition}

\subsection{Aubert dual of tempered packets}\label{sec:subregular-aubert}

\begin{proposition}\label{prop:subregular-abelian-Aubert-dual-computation}
      In $R(\rG_2)$, we have the following results on the Aubert duals of the tempered subregular $L$-packets of ${\rG_2}$.

      \begin{enumerate}
            \item We have
            \begin{align*}
                  \widehat{\pi_{\gen}[1]} &= J_Q(\pi(1, 1) \nu_2), \\
                  \widehat{\pi_{\deg}[1]} &= J_P(\St_{\GL_2}\nu_2^{1/2}), \\
                  \widehat{J_Q(\St_{\GL_2}\nu_2^{1/2})} &= J_Q(\St_{\GL_2}\nu_2^{1/2}).
            \end{align*}
            Furthermore,
            \begin{align*}
                  [I_P(\mathbf{1}_{\GL_2}\nu_2^{1/2})] = [\pi_{\deg}[1]] + [\widehat{\pi_{\gen}[1]}] + [J_Q(\St_{\GL_2}\nu_2^{1/2})], \\
                  [I_P(\St_{\GL_2}\nu_2^{1/2})]
                  = [\pi_{\gen}[1]] + [\widehat{\pi_{\deg}[1]}] + [J_Q(\St_{\GL_2}\nu_2^{1/2})], \\
                  [I_Q(\St_{\GL_2}\nu_2^{1/2})] = [\pi_{\gen}[1]] + [\pi_{\deg}[1]] + [J_Q(\St_{\GL_2}\nu_2^{1/2})],
            \end{align*}
            and
            \begin{align*}
                  [J_Q(\pi(1, 1) \nu_2)] &= [I_Q(\pi(1, 1)\nu_2)] - [I_P(\St_{\GL_2}\nu_2^{1/2})] - [I_Q(\St_{\GL_2}\nu_2^{1/2})] + [\pi_{\mathrm{gen}}[1]], \\
                  [J_P(\St_{\GL_2}\nu_2^{1/2})] &= [I_P(\St_{\GL_2}\nu_2^{1/2})] - [I_Q(\St_{\GL_2}\nu_2^{1/2})] + [\pi_{\deg}[1]]. \\
            \end{align*}

            \item Let $\chi$ be a non-trivial quadratic character of $F^{\times}$.
            We have
            \begin{align*}
                  \widehat{\pi_{\gen}[\chi]} &= J_Q(\pi(1, \chi) \nu_2), \\
                  \widehat{J_P(\St_2(\chi)\nu_2^{1/2})} &= J_Q(\St_2(\chi)\nu_2^{1/2}).
            \end{align*}
            Furthermore,
            \begin{align*}
                  [I_P(\chi \circ \det \nu_2^{1/2})] = [\widehat{\pi_{\gen}[\chi]}] + [J_Q(\St_2(\chi)\nu_2^{1/2})], \\
                  [I_P(\St_2(\chi)\nu_2^{1/2})] = [\pi_{\gen}[\chi]] + [J_P(\St_2(\chi)\nu_2^{1/2})], \\
                  [I_Q(\St_2(\chi)\nu_2^{1/2})] = [\pi_{\gen}[\chi]] + [J_Q(\St_2(\chi)\nu_2^{1/2})].
            \end{align*}
            and
            \begin{align*}
                  [J_Q(\pi(1, \chi) \nu_2)] &= [I_Q(\pi(1, \chi)\nu_2)] - [I_P(\St_2(\chi)\nu_2^{1/2})] - [I_Q(\St_{2}(\chi)\nu_2^{1/2})] + [\pi_{\mathrm{gen}}[\chi]].
            \end{align*}

            \item Let $\chi$ be a non-trivial cubic character of $F^{\times}$.
            We have
            \begin{align*}
                  \widehat{\pi_{\gen}[\chi]} &= J_Q(\pi(\chi, \chi^{-1}) \nu_2), \\
                  \widehat{J_P(\St_2(\chi)\nu_2^{1/2})} &= J_P(\St_2(\chi^{-1})\nu_2^{1/2}).
            \end{align*}
            Furthermore,
            \begin{align*}
                  [I_P(\chi \circ \det \nu_2^{1/2})] &=
                  [\widehat{\pi_{\gen}[\chi]}] + [J_P(\St_2(\chi^{-1})\nu_2^{1/2})], \\
                  [I_P(\St_2(\chi^{-1})\nu_2^{1/2})] &=
                  [\pi_{\gen}[\chi]] + [J_P(\St_2(\chi^{-1})\nu_2^{1/2})],
            \end{align*}
            and we have
            \begin{align*}
                  [J_Q(\pi(\chi, \chi^{-1}) \nu_2)]
                   &= [I_Q(\pi(\chi, \chi^{-1})\nu_2)] - [I_P(\St_2(\chi)\nu_2^{1/2})] \\
                   &- [I_P(\St_{2}(\chi^{-1})\nu_2^{1/2})] + [\pi_{\mathrm{gen}}[\chi]].
            \end{align*}
      \end{enumerate}
\end{proposition}

\begin{proof}
      The composition series are due to Mui\'c \cite{Muic1997-UnitaryDualG2}*{Propositions 4.3, 4.1, and 4.2}.
      We use their formulations in \cite{Qin2026-AubertZelevinskyInvolutionG2HeckeAlgebras}*{ Propositions 6.2.1(2), 7.4.1(2), and 7.5.1(2)}, with the inducing data expressed in the intrinsic $\GL_2$ coordinates of Section~\ref{subsec:Qin-notation}.
      In the first case Qin's $\pi'(1)$ is $\pi_{\gen}[1]$ and $\pi(1)$ is $\pi_{\deg}[1]$; the identification also follows from the generic and nongeneric subrepresentations in \cite{GanSavin2023-HoweDualityDichotomyExceptionalThetaCorrespondences}*{ Propositions 3.1 and 3.2}.
      The corresponding duality assertions are \cite{Qin2026-AubertZelevinskyInvolutionG2HeckeAlgebras}*{ Propositions 6.2.3, 7.4.3, and 7.5.3}.
      Qin uses $\sum_{I\subset\Delta}(-1)^{|I|}I_{P_I}r_{P_I}$, which equals \eqref{eq:Aubert-operator} since $\operatorname{rank}{\rG_2}=2$.
      All these principal-series constituents have toral support, so their support sign is $+1$.
      Thus these duality assertions agree with positive Aubert duality.
      The supercuspidal case follows from Proposition~\ref{prop:Aubert}.
\end{proof}

\begin{proposition}
      Assume that $\sigma(W_F) = S_3$.
      Let $\tau$ be a supercuspidal representation of $\GL_2(F)$ corresponding to the $L$-parameter $\rho \circ \sigma \colon W_F \to \GL_2(\C)$.
      We have an exact sequence
      \begin{align*}
            0 \to \pi_{\gen}[\tau] \to I_Q(\tau \nu_2) \to J_Q(\tau \nu_2) \to 0.
      \end{align*}

      We also have
      \begin{align*}
            \widehat{\pi_{\gen}[\tau]} &= J_Q(\tau \nu_2).
      \end{align*}
\end{proposition}

\begin{proof}
      This follows from \cite{GanSavin2023-HoweDualityDichotomyExceptionalThetaCorrespondences}*{Proposition 3.2 (i), the second bullet point}.
\end{proof}

\subsection{The case with trivial image}\label{sec:subregular-trivial}

We consider the case where $\sigma(W_F) = \{1\}$.
In this case, we have $\cS_{\psi} = S_3$.
Let $\psi'$ (resp. $\phi'$) be the lift of $\psi$ (resp. $\phi$) to $\PGSO_8(F)$.

\begin{proposition}\label{prop:subregular-trivial-A-packet-lift-PGSO8}
      The packet $\Pi_{\psi'}$ consists of two representations of $\PGSO_8(F)$
      \begin{align*}
            J_{Q'}((\pi(1, 1)\nu_2)_{M_{Q'}})
            \quad
            \text{and}
            \quad
            J_{P'} ((\St_{\GL_2}\nu_2^{1/2})_{M_{P'}}).
      \end{align*}
\end{proposition}

\begin{proof}
      This follows from \cite{BanJantzen2003-DegeneratePrincipalSeriesEvenOrthogonalGroups}*{Theorem 5.3, (3)(b)} with $k=3, l=1, n=0, \rho = \mathbf{1}, \sigma = \mathbf{1}$ and $\alpha = 0$ and Proposition \ref{prop:Xu-Aubert} since the cuspidal support of these representations coincides with the cuspidal support of the lift of $\phi$.
\end{proof}

\begin{proposition}\label{prop:subregular-trivial-ecr}
      The distribution
      \begin{align*}
            S\Theta_{\psi} := \Theta_{J_Q(\pi(1, 1) \nu_2)} + 2\Theta_{J_P(\St_{\GL_2}\nu_2^{1/2})} + \Theta_{\pi_{\mathrm{sc}}[1]}
      \end{align*}
      is equal to
      \begin{align*}
            \Theta_{I_Q(\pi(1, 1) \nu_2)} + \Theta_{I_P(\St_{\GL_2}\nu_2^{1/2})} -3 \Theta_{I_Q(\St_{\GL_2}\nu_2^{1/2})} + S\Theta_{\phi}
      \end{align*}
      and is stable.
      This transfers to the distribution
      \begin{align*}
            \cT S\Theta_{\psi} =
            \Theta_{\widetilde{J_{Q'}((\pi(1, 1)\nu_2)_{M_{Q'}})}}
            + \Theta_{\widetilde{J_{P'} ((\St_{\GL_2}\nu_2^{1/2})_{M_{P'}})}}.
      \end{align*}

      We set
      \begin{align*}
            \pi_{\psi}(\mathbf{1}) &= J_Q(\pi(1, 1) \nu_2) \\
            \pi_{\psi}(\rho) &= J_P(\St_{\GL_2}\nu_2^{1/2}) \\
            \pi_{\psi}(\epsilon) &= \pi_{\mathrm{sc}}[1].
      \end{align*}
      Then, we have
      \begin{align*}
            \Trans^{\rG_2, *}_{\SO_4}(\Theta_{\mathbf{1}_{\SO_4}}) &= \Theta_{\pi_{\psi}(\mathbf{1})} - \Theta_{\pi_{\psi}(\epsilon)}, \\
            \Trans^{\rG_2, *}_{\PGL_3}(\Theta_{\mathbf{1}_{\PGL_3}}) &= \Theta_{\pi_{\psi}(\mathbf{1})} - \Theta_{\pi_{\psi}(\rho)} + \Theta_{\pi_{\psi}(\epsilon)}.
      \end{align*}
\end{proposition}

\begin{proof}
      This follows from Proposition \ref{prop:subregular-abelian-Aubert-dual-computation} (1), Proposition \ref{prop:subregular-trivial-A-packet-lift-PGSO8} and Theorem \ref{thm:signed-transfer}.
\end{proof}

\subsection{The cases with quadratic and cubic image}\label{sec:subregular-quadratic}

We now consider the case where $\sigma(W_F) = \mu_2$ or $\mu_3$.

\begin{proposition}\label{prop:subregular-quadratic-A-packet-lift-SO8}
      If $\chi$ is a non-trivial quadratic character of $F^{\times}$, then the packet $\Pi_{\phi'}$ of $\PGSO_8(F)$ consists of discrete series representations $\pi'_{\gen}[\chi]$ and $\pi'_{\deg}[\chi]$ corresponding to the two irreducible representations $\mathbf{1}$ and $\epsilon$ of $\mu_2$.
      The latter is supercuspidal.
      Thus, the packet $\Pi_{\psi'}$ of $\PGSO_8(F)$ consists of
      \begin{align*}
            J_{Q'}
            ((\pi(1, \chi)\nu_2)_{M_{Q'}}) \quad \text{and} \quad \pi'_{\deg}[\chi].
      \end{align*}
\end{proposition}

\begin{proof}
      This follows from Lemma \ref{lem:triality-invariant-non-tempered}, \cite{Xu2017-MoeglinParametrizationArthurPackets}*{Definition 6.3(i), (iii)(c)(1), p.920} and \cite{Xu2017-CuspidalSupportDiscreteSeries}*{Theorem 3.3}.
\end{proof}

\begin{proposition}\label{prop:subregular-cubic-A-packet-lift-SO8}
      Let $\chi$ be a non-trivial cubic character of $F^{\times}$.
      Then, the packet $\Pi_{\phi'}$ of $\PGSO_8(F)$ is a singleton consisting of the representation
      \begin{align*}
            I^{\PGSO_8}_{\GL_3 \times \GL_1 \times \GL_1/\GL_1}(\St_{3}(\chi) \boxtimes \mathbf{1} \boxtimes \mathbf{1}).
      \end{align*}
      Thus, the packet $\Pi_{\psi'}$ of $\PGSO_8(F)$ is also a singleton consisting of the representation
      \begin{align*}
            I^{\PGSO_8}_{\GL_3 \times \GSO_2 / \GL_1}(\chi \circ \det \boxtimes \mathbf{1}_{\GSO_2})
            = I^{\PGSO_8}_{\GL_3 \times \GSO_2 / \GL_1}(\chi^{-1} \circ \det \boxtimes \mathbf{1}_{\GSO_2}).
      \end{align*}
      This representation is also equal to
      \begin{align*}
            J^{\PGSO_8}_{Q'}(((\pi(\chi, \chi^{-1})\nu_2)_{M_{Q'}})).
      \end{align*}
\end{proposition}

\begin{proof}
      The first assertion follows from the fact that the standard parameter of $\phi'$ is given as
      \begin{align*}
            \chi \boxtimes S_3^D \oplus \chi^{-1} \boxtimes S_3^D \oplus 2 \cdot \mathbf{1}
      \end{align*}
      The second assertion follows from Proposition \ref{prop:Aubert}.
      The last assertion follows from \cite{Xu2024-ArthurPacketsQuasisplitGSp2nGO2n}*{Lemma 9.3} and Lemma \ref{lem:triality-invariant-non-tempered}.
\end{proof}

Then, we obtain the following results from Theorem \ref{thm:signed-transfer} as in the previous subsection for the case with trivial image.

\begin{proposition}\label{prop:subregular-quadratic-ecr}
      For a non-trivial quadratic character $\chi$, the distribution
      \begin{align*}
            S\Theta_{\psi} = \Theta_{J_Q(\pi(1, \chi) \nu_2)} + \Theta_{\pi_{\mathrm{sc}}[\chi]}
      \end{align*}
      is stable and equal to
      \begin{align*}
            \Theta_{I_Q(\pi(1, \chi) \nu_2)} - \Theta_{I_P(\St_2(\chi)\nu_2^{1/2})} - \Theta_{I_Q(\St_2(\chi)\nu_2^{1/2})} + S\Theta_{\phi}.
      \end{align*}
      This distribution transfers to the distribution
      \begin{align*}
            \cT S\Theta_{\psi}
            = \Theta_{\widetilde{
                  J_{Q'}
            ((\pi(1, \chi)\nu_2)_{M_{Q'}})
            }}
            + \Theta_{\widetilde{\pi'_{\deg}[\chi]}}.
      \end{align*}
      Furthermore, we have
      \begin{align*}
            \Theta_{\widetilde{
             J^{\PGSO_8}_{Q'}(((\pi(1, \chi)\nu_2)_{M_{Q'}}))
            }}
            - \Theta_{\widetilde{\pi'_{\deg}[\chi]}}
            = \Trans^{\widetilde{\PGSO_8}, *}_{\SO_4}(\Theta_{\chi \circ \sn}).
      \end{align*}

      We set
      \begin{align*}
            \pi_{\psi}(\mathbf{1}) &= J_Q(\pi(1, \chi) \nu_2) \\
            \pi_{\psi}(\epsilon) &= \pi_{\mathrm{sc}}[\chi].
      \end{align*}
      Then, we have
      \begin{align*}
            \Trans^{\rG_2, *}_{\SO_4}(\Theta_{\chi \circ \sn}) &=
            \Theta_{\pi_{\psi}(\mathbf{1})} - \Theta_{\pi_{\psi}(\epsilon)}.
      \end{align*}
\end{proposition}

\begin{proposition}\label{prop:subregular-cubic-ecr}
      For a non-trivial cubic character $\chi$, the distribution
      \begin{align*}
            S\Theta_{\psi} = \Theta_{J_Q(\pi(\chi, \chi^{-1}) \nu_2)} + \Theta_{\pi_{\mathrm{sc}, 1}[\chi]} + \Theta_{\pi_{\mathrm{sc}, 2}[\chi]}
      \end{align*}
      is stable and equal to
      \begin{align*}
            \Theta_{I_Q(\pi(\chi, \chi^{-1}) \nu_2)} - \Theta_{I_P(\St_2(\chi)\nu_2^{1/2})} - \Theta_{I_P(\St_2(\chi^{-1})\nu_2^{1/2})} + S\Theta_{\phi}.
      \end{align*}
      This distribution transfers to the distribution
      \begin{align*}
            \cT S\Theta_{\psi}
            = \Theta_{\widetilde{
                 J^{\PGSO_8}_{Q'}(((\pi(\chi, \chi^{-1})\nu_2)_{M_{Q'}}))
            }}.
      \end{align*}

      We set
      \begin{align*}
            \pi_{\psi}(\mathbf{1}) &= J_Q(\pi(\chi, \chi^{-1}) \nu_2) \\
            \pi_{\psi}(\chi) &= \pi_{\mathrm{sc}, 1}[\chi] \\
            \pi_{\psi}(\chi^2) &= \pi_{\mathrm{sc}, 2}[\chi].
      \end{align*}
      Then, we have
      \begin{align*}
            \Trans^{\rG_2, *}_{\PGL_3}(\Theta_{\chi \circ \det}) &= \Theta_{\pi_{\psi}(\mathbf{1})} + \zeta \Theta_{\pi_{\psi}(\chi)} + \zeta^2 \Theta_{\pi_{\psi}(\chi^2)}, \\
            \Trans^{\rG_2, *}_{\PGL_3}(\Theta_{\chi^2 \circ \det}) &= \Theta_{\pi_{\psi}(\mathbf{1})} + \zeta^2 \Theta_{\pi_{\psi}(\chi)} + \zeta \Theta_{\pi_{\psi}(\chi^2)}.
      \end{align*}
\end{proposition}

\subsection{The case with $S_3$-image}\label{sec:subregular-full}

Let $\tau$ be the supercuspidal representation of $\GL_2(F)$ corresponding to the $L$-parameter $\rho \circ \sigma \colon W_F \to \GL_2(\C)$.
Let $\omega$ be the central character of $\tau$.

Then, we have the following result.
\begin{proposition}\label{prop:subregular-full-ecr}
      The packet $\Pi_{\psi'}$ of $\PGSO_8(F)$ is a singleton consisting of the representation
      \begin{align*}
            J_{Q'}((\tau \nu_2)_{M_{Q'}}).
      \end{align*}

      The distribution
      \begin{align*}
            \Theta_{J_{Q}(\tau \nu_2)}
      \end{align*}
      is stable and transfers to the distribution
      \begin{align*}
            \cT \Theta_{J_{Q}(\tau \nu_2)} = \Theta_{\widetilde{J_{Q'}((\tau \nu_2)_{M_{Q'}})}}.
      \end{align*}

      We set
      \begin{align*}
            \pi_{\psi}(\mathbf{1}) = J_{Q}(\tau \nu_2).
      \end{align*}
\end{proposition}

\begin{proof}
      This follows from Theorem \ref{thm:signed-transfer} with the support sign $-1$ on both sides.
\end{proof}

\section{Co-tempered long-root and short-root parameters: supercuspidal cases}\label{sec:co-tempered-long-short-supercuspidal}

We now consider the co-tempered long-root and short-root parameters of ${\rG_2}$.
This is given by a supercuspidal irreducible representation $\sigma$ of $\PGL_2(F)$.
The swapped parameter $\phi = \widehat\psi$ is tempered with $\cS_{\phi} \cong \mu_2$.
Thus, the $L$-packet $\Pi_{\phi}$ consists of a generic representation $\delta = \delta_{R}(\sigma)$ in Lemma \ref{lem:rank-one} for $R \in \{P, Q\}$ and $s = 1/2$ (see also \cite{GanSavin2023-HoweDualityDichotomyExceptionalThetaCorrespondences}*{Propositions 3.1 and 3.2}) and a non-generic supercuspidal representation $\pi_{\mathrm{sc}}[\sigma]$ by \cite{Xu2017-CuspidalSupportDiscreteSeries}*{Theorem 3.3} and \cite{GanSavin2023-LocalLanglandsConjectureG2}*{Proposition 8.2}.

Furthermore, the $L$-packet attached to the lift of $\phi$ to $\PGSO_8(F)$ consists of two discrete series representations, namely a generic representation $\delta'$ and a non-generic supercuspidal representation $\pi'_{\mathrm{sc}}[\sigma]$.
These representations are obtained from the composition of the exceptional theta lift from $\rG_2$ to $\PGSp_6$ and the similitude theta lift from $\PGSp_6$ to $\PGSO_8$.

We can compute the Aubert dual of $\delta$ in Lemma \ref{lem:rank-one} and obtain the endoscopic character identities by Theorem \ref{thm:signed-transfer}.

\begin{proposition}\label{prop:co-tempered-long-short-ecr}
      We set
      \begin{align*}
            \pi_{\psi}(\mathbf{1}) &= \widehat{\delta}, \\
            \pi_{\psi}(\epsilon) &= \pi_{\mathrm{sc}}[\sigma].
      \end{align*}
      Then, the distribution
      \begin{align*}
            S\Theta_{\psi} = \Theta_{\pi_{\psi}(\mathbf{1})} - \Theta_{\pi_{\psi}(\epsilon)}
      \end{align*}
      is stable and transfers to the distribution
      \begin{align*}
            \cT S\Theta_{\psi}
            = \Theta_{\widetilde{\widehat{\delta'}}} - \Theta_{\widetilde{\pi'_{\mathrm{sc}}[\sigma]}}.
      \end{align*}

      If $\psi$ is a long root (resp. short root) parameter, then we have
      \begin{align*}
            \delta = \delta_{Q}(\sigma) \quad (\text{resp. } \delta = \delta_{P}(\sigma)),
      \end{align*}
      in \cite{GanSavin2023-HoweDualityDichotomyExceptionalThetaCorrespondences}.

      Furthermore, we have
      \begin{align*}
            \Theta_{\widetilde{\widehat{\delta'}}} + \Theta_{\widetilde{\pi'_{\mathrm{sc}}[\sigma]}}
            = \Trans^{\widetilde{\PGSO_8}, *}_{\SO_4}(\Theta_{J^{\SO_4}_{R}(\sigma \nu_2^{1/2})}), \\
            \Theta_{\pi_{\psi}(\mathbf{1})} + \Theta_{\pi_{\psi}(\epsilon)}
            = \Trans^{\rG_2, *}_{\SO_4}(\Theta_{J^{\SO_4}_{R}(\sigma \nu_2^{1/2})}),
      \end{align*}
      where $R = P_{l}$ (resp. $R = P_{s}$) if $\psi$ is a long root (resp. short root) parameter and $J^{\SO_4}_{R}(\sigma \nu_2^{1/2})$ is the Langlands quotient.
\end{proposition}

\section{A-parameters with unbounded image}\label{sec:A-parameter-unbounded-ecr}

In the long-root and short-root cases, we can consider the A-parameters with unbounded image associated with generic unitary principal series representations of $\PGL_2(F)$.
Let $\sigma$ be a generic unitary principal series irreducible representation of $\PGL_2(F)$.
We consider the long-root and short-root A-parameter $\psi$ of $\rG_2(F)$ associated with $\sigma$.
The long-root parameter is of the form
\begin{align*}
      \psi \colon W_F \times \SL_2^A(\C) \xrightarrow{\phi_{\sigma} \times \id_{\SL_2}} \SL_{2, s}(\C) \times \SL_{2, l}(\C) \to \rG_2(\C).
\end{align*}
For the short-root parameter, interchange the subscripts $s$ and $l$ in this formula.
Here, we relax the boundedness condition on the image of $\phi_{\sigma}$.

We first consider the long-root case.
Then, the parameter factors through the long root Levi subgroup of $\rG_2(\C)$ and defines an A-parameter for $\GL_2(F)$.
This corresponds to a character $\chi_{M_P}$ of $M_P$.
In this case, we set
\begin{align*}
      \pi_{\psi}(\mathbf{1}) = I_P(\chi_{M_P}).
\end{align*}
By \cite{Muic1997-UnitaryDualG2}*{Theorem 3.1 (i), Lemma 5.2} the representation $I_P(\chi_{M_P})$ is irreducible and unitary.
We have the twisted endoscopic character identity
\begin{align*}
      \Trans^{\widetilde{\PGSO_8}, *}_{\rG_2} \Theta_{I_P(\chi_{M_P})} = \Theta_{\widetilde{I_{P'}(\chi_{M_{P'}})}},
\end{align*}
where $\chi_{M_{P'}}$ is the character of $M_{P'}(F)$ corresponding to $\chi_{M_P}$.

The short-root case is similar: The parameter factors through the short root Levi subgroup of $\rG_2(\C)$ and defines an A-parameter for $\GL_2(F)$.
This corresponds to a character $\chi_{M_Q}$ of $M_Q$.
In this case, we set
\begin{align*}
      \pi_{\psi}(\mathbf{1}) = I_Q(\chi_{M_Q}).
\end{align*}
By \cite{Muic1997-UnitaryDualG2}*{Theorem 3.1 (ii), Lemma 5.2} the representation $I_Q(\chi_{M_Q})$ is irreducible and unitary.
We have the twisted endoscopic character identity
\begin{align*}
      \Trans^{\widetilde{\PGSO_8}, *}_{\rG_2} \Theta_{I_Q(\chi_{M_Q})} = \Theta_{\widetilde{I_{Q'}(\chi_{M_{Q'}})}},
\end{align*}
where $\chi_{M_{Q'}}$ is the character of $M_{Q'}(F)$ corresponding to $\chi_{M_Q}$.

We summarize the above discussion in the following proposition.
\begin{proposition}\label{prop:co-tempered-long-short-unbounded-ecr}
      Let $\psi$ be a long-root (resp. short-root) A-parameter of $\rG_2(F)$ associated with a generic unitary principal series representation $\sigma$ of $\PGL_2(F)$.
      Then, the A-packet $\Pi_{\psi}$ consists of a single representation $\pi_{\psi}(\mathbf{1})$ given by $I_P(\chi_{M_P})$ (resp. $I_Q(\chi_{M_Q})$).
      The distribution $\Theta_{\pi_{\psi}(\mathbf{1})}$ is stable and transfers to the distribution $\Theta_{\widetilde{I_{P'}(\chi_{M_{P'}})}}$ (resp. $\Theta_{\widetilde{I_{Q'}(\chi_{M_{Q'}})}}$).
\end{proposition}

\section{Mixed long root parameters}\label{sec:mixed-long}

Let $\psi_D$ be a mixed long root parameter of ${\rG_2}$ associated with a quadratic character $\chi$ of $F^{\times}$.
Let $\psi'_D$ denote the lift to $\PGSO_8(F)$.

\begin{lemma}\label{lem:mixed-long-lift}
      There is a triality-invariant A-packet corresponding to $\psi'_D$ for $\PGSO_8(F)$ consisting of two representations.
      One is non-tempered and the other is tempered.

      If $\chi$ is trivial, then the non-tempered representation is given by
      \begin{align*}
            J_{Q'}((\St_{\GL_2} \nu_2^{1/2})_{M_{Q'}}),
      \end{align*}
      and the tempered representation has the standard parameter
      \begin{align*}
            2 \cdot S_3^D \oplus 2 \cdot \mathbf{1}
      \end{align*}
      and it is non-generic.

      If $\chi$ is non-trivial, then the non-tempered representation is given by
      \begin{align*}
            J_{Q'}((\St_{2}(\chi) \nu_2^{1/2})_{M_{Q'}}),
      \end{align*}
      and the other representation is supercuspidal and non-generic with the standard parameter
      \begin{align*}
            \chi \boxtimes S_3^D \oplus \chi \oplus S_3^D \oplus \mathbf{1}.
      \end{align*}
\end{lemma}

\begin{proof}
      We first determine the A-packet for $\SO_8(F)$ corresponding to the lift of $\psi_D$.
      If $\chi$ is trivial, then \cite{Xu2021-NonarchimedeanComponentsNonEndoscopicAutomorphicRepresentations}*{Theorem 1.3 (2)} implies the A-packet consists of two representations, one is non-tempered and given by
      \begin{align*}
            J_{\GL_2 \times \SO_4}(\St_{\GL_2} \nu_2^{1/2} \boxtimes \St_{\SO_4})
      \end{align*}
      and the other is a tempered and non-generic representation $\pi'_{\temp, \mathrm{ng}}$ with the standard parameter $2 \cdot S_3^D \oplus 2 \cdot \mathbf{1}$.
      If $\chi$ is non-trivial, then we can use \cite{Xu2017-MoeglinParametrizationArthurPackets}*{Theorem 7.8}, which implies that the A-packet for $\SO_8(F)$ consists of three representations, one is non-tempered and given by
      \begin{align*}
            J_{\GL_2 \times \SO_4}(\St_{2}(\chi) \nu_2^{1/2} \boxtimes \St_{\SO_4})
      \end{align*}
      and the others are supercuspidal and non-generic with the standard parameter $\chi \boxtimes S_3^D \oplus \chi \oplus S_3^D \oplus \mathbf{1}$.
      By \cite{Xu2018-LPacketsQuasisplitGSp2nGO2n}*{Lemma 6.6}, these supercuspidal representations form a single $\PGSO_8(F)$-orbit $\pi'_{\mathrm{sc}, \chi}$.

      These observations and Lemma \ref{lem:triality-invariant-non-tempered} imply that if there is an A-packet with the standard parameter $\Std_{\PGSO_8} \circ \psi'_D$ and triality-invariant, then it is the A-packet for $\PGSO_8(F)$ containing
      \begin{align*}
            J_{Q'}((\St_{\GL_2} \nu_2^{1/2})_{M_{Q'}})
      \end{align*}
      if $\chi$ is trivial and
      \begin{align*}
            J_{Q'}((\St_{2}(\chi) \nu_2^{1/2})_{M_{Q'}})
      \end{align*}
      if $\chi$ is non-trivial.
      We now show that they actually form an A-packet for $\PGSO_8(F)$.
      By \cite{Xu2024-ArthurPacketsQuasisplitGSp2nGO2n}*{Theorem 8.1 (1)}, it suffices to show the stability of the distribution
      \begin{align*}
            \Theta_{J_{Q'}((\St_{\GL_2} \nu_2^{1/2})_{M_{Q'}})} - \Theta_{\pi'_{\temp, \mathrm{ng}}}, \\
            \Theta_{J_{Q'}((\St_{2}(\chi) \nu_2^{1/2})_{M_{Q'}})} - \Theta_{\pi'_{\mathrm{sc}, \chi}}.
      \end{align*}
      For $\chi = \mathbf{1}$, we can apply \cite{Muic2005-ReducibilityGeneralizedPrincipalSeries}*{Lemma 5.2} to see that the parabolic induction $I_{Q'}((\St_{\GL_2} \nu_2^{1/2})_{M_{Q'}})$ is of length two and has a unique irreducible tempered generic subrepresentation $\pi'_{\temp, \gen}$.
      The standard parameter of this representation must be $2 \cdot S_3^D \oplus 2 \cdot \mathbf{1}$.
      Thus, we have
      \begin{align*}
            [J_{Q'}((\St_{\GL_2} \nu_2^{1/2})_{M_{Q'}})] - [\pi'_{\temp, \mathrm{ng}}]
            = [I_{Q'}((\St_{\GL_2} \nu_2^{1/2})_{M_{Q'}})] - [\pi'_{\temp, \gen}] - [\pi'_{\temp, \mathrm{ng}}].
      \end{align*}
      The corresponding distribution is stable by Theorem \ref{thm:tempered-input}.

      The case for non-trivial $\chi$ is similar by using \cite{Muic2005-ReducibilityStandardRepresentations}*{Lemma 6.1}.
      In this case, we obtain the exact sequence
      \begin{align*}
            0 \to \pi'_{\gen, \chi} \to I_{Q'}((\St_{2}(\chi) \nu_2^{1/2})_{M_{Q'}}) \to J_{Q'}((\St_{2}(\chi) \nu_2^{1/2})_{M_{Q'}}) \to 0,
      \end{align*}
      where $\pi'_{\gen, \chi}$ is the unique irreducible generic discrete series subrepresentation of $I_{Q'}((\St_{2}(\chi) \nu_2^{1/2})_{M_{Q'}})$ with the standard parameter $\chi \boxtimes S_3^D \oplus \chi \oplus S_3^D \oplus \mathbf{1}$.
      Thus, we have
      \begin{align*}
            [J_{Q'}((\St_{2}(\chi) \nu_2^{1/2})_{M_{Q'}})] - [\pi'_{\mathrm{sc}, \chi}]
            = [I_{Q'}((\St_{2}(\chi) \nu_2^{1/2})_{M_{Q'}})] - ([\pi'_{\gen, \chi}] + [\pi'_{\mathrm{sc}, \chi}]).
      \end{align*}
      The corresponding distribution is stable by Theorem \ref{thm:tempered-input}.
\end{proof}

We record the following corollary of the proof above.
\begin{corollary}\label{cor:mixed-long-exact}
      We have the following exact sequences
      \begin{align*}
            0 \to \pi'_{\temp, \gen} \to I_{Q'}((\St_{\GL_2} \nu_2^{1/2})_{M_{Q'}}) \to J_{Q'}((\St_{\GL_2} \nu_2^{1/2})_{M_{Q'}}) \to 0, \\
            0 \to \pi'_{\mathrm{gen}, \chi} \to I_{Q'}((\St_{2}(\chi) \nu_2^{1/2})_{M_{Q'}}) \to J_{Q'}((\St_{2}(\chi) \nu_2^{1/2})_{M_{Q'}}) \to 0.
      \end{align*}
\end{corollary}

\Needspace{10\baselineskip}
\begin{proposition}\label{prop:mixed-long-root-ecr}
      Let $\psi_D$ be a mixed long root parameter of ${\rG_2}$ associated with a quadratic character $\chi$ of $F^{\times}$.
      \begin{enumerate}
            \item If $\chi$ is trivial, then the distribution
            \begin{align*}
                  S\Theta_{\psi_D} = \Theta_{J_Q(\St_{\GL_2} \nu_2^{1/2})} - (\Theta_{\pi_{\deg}[1]} + \Theta_{\pi_{\mathrm{sc}}[1]}),
            \end{align*}
            is stable and it is equal to
            \begin{align*}
                  \Theta_{I_Q(\St_{\GL_2} \nu_2^{1/2})} - (\Theta_{\pi_{\gen}[1]} + 2\Theta_{\pi_{\deg}[1]} + \Theta_{\pi_{\mathrm{sc}}[1]}).
            \end{align*}
            This distribution transfers to the distribution
            \begin{align*}
                  \cT S\Theta_{\psi_D} =
                  \Theta_{\widetilde{J_{Q'}((\St_{\GL_2} \nu_2^{1/2})_{M_{Q'}})}} - \Theta_{\widetilde{\pi'_{\temp, \mathrm{ng}}}},
            \end{align*}
            where $\pi'_{\temp, \mathrm{ng}}$ is the triality-invariant tempered non-generic representation of $\PGSO_8(F)$ with the standard parameter $2 \cdot S_3^D \oplus 2 \cdot \mathbf{1}$.

            Furthermore, we have distributions
            \begin{align*}
                  \Theta_{J_Q(\St_{\GL_2} \nu_2^{1/2})} + (\Theta_{\pi_{\deg}[1]} + \Theta_{\pi_{\mathrm{sc}}[1]})
            \end{align*}
            and
            \begin{align*}
                  \Theta_{\widetilde{J_{Q'}((\St_{\GL_2} \nu_2^{1/2})_{M_{Q'}})}} + \Theta_{\widetilde{\pi'_{\temp, \mathrm{ng}}}},
            \end{align*}
            which are transfers of the stable distribution
            \begin{align*}
                  \Theta_{J_{P_{l}}^{\SO_4}(\St_{\GL_2} \nu_2^{1/2})}
            \end{align*}
            of $\SO_4(F)$.
            Here, we recall that $P_{l}$ is the long root parabolic subgroup of $\SO_4(F)$ with the short root dual parabolic subgroup.

            We set
            \begin{align*}
                  \pi_{\psi_D}(\mathbf{1}) &= J_Q(\St_{\GL_2} \nu_2^{1/2}), \\
                  \pi_{\psi_D}(\epsilon) &= \pi_{\deg}[1] \oplus \pi_{\mathrm{sc}}[1].
            \end{align*}

            \item If $\chi$ is non-trivial, then the distribution
            \begin{align*}
                  S\Theta_{\psi_D} =
                  \Theta_{J_Q(\St_{2}(\chi) \nu_2^{1/2})} - \Theta_{\pi_{\mathrm{sc}}[\chi]},
            \end{align*}
            is stable and it is equal to
            \begin{align*}
                  \Theta_{I_Q(\St_{2}(\chi) \nu_2^{1/2})} - (\Theta_{\pi_{\gen}[\chi]} + \Theta_{\pi_{\mathrm{sc}}[\chi]}).
            \end{align*}
            This distribution transfers to the distribution
            \begin{align*}
                  \cT S\Theta_{\psi_D} =
                  \Theta_{\widetilde{J_{Q'}((\St_{2}(\chi) \nu_2^{1/2})_{M_{Q'}})}} - \Theta_{\widetilde{\pi'_{\mathrm{sc}, \chi}}},
            \end{align*}
            where $\pi'_{\mathrm{sc}, \chi}$ is the triality-invariant supercuspidal non-generic representation of $\PGSO_8(F)$ with the standard parameter
            \begin{align*}
                  \chi \boxtimes S_3^D \oplus \chi \oplus S_3^D \oplus \mathbf{1}.
            \end{align*}

            Furthermore, we have distributions
            \begin{align*}
                  \Theta_{J_Q(\St_{2}(\chi) \nu_2^{1/2})} + \Theta_{\pi_{\mathrm{sc}}[\chi]}
            \end{align*}
            and
            \begin{align*}
                  \Theta_{\widetilde{J_{Q'}((\St_{2}(\chi) \nu_2^{1/2})_{M_{Q'}})}} + \Theta_{\widetilde{\pi'_{\mathrm{sc}, \chi}}},
            \end{align*}
            which are transfers of the stable distribution
            \begin{align*}
                  \Theta_{J_{P_{l}}^{\SO_4}(\St_{2}(\chi) \nu_2^{1/2})}
            \end{align*}
            of $\SO_4(F)$.

            We set
            \begin{align*}
                  \pi_{\psi_D}(\mathbf{1}) &= J_Q(\St_{2}(\chi) \nu_2^{1/2}), \\
                  \pi_{\psi_D}(\epsilon) &= \pi_{\mathrm{sc}}[\chi].
            \end{align*}
      \end{enumerate}
\end{proposition}

\begin{proof}
      The stability follows from Theorem \ref{thm:tempered-input} and the transfer follows from Theorem \ref{thm:signed-transfer} and Theorem \ref{thm:tempered-input}.
      In these cases, the component group $\cS_{\psi'_D}$ is of order two and the corresponding endoscopic group is $\SO_4(F)$.
\end{proof}

\section{Mixed short root parameters}\label{sec:mixed-short}

By applying the (twisted) Aubert duality, we can obtain the results for the mixed short root parameters of ${\rG_2}$ from the results for the mixed long root parameters in Proposition \ref{prop:mixed-long-root-ecr}.

\begin{proposition}\label{prop:mixed-short-root-trivial-ecr}
      We have the following results for the mixed short root parameters of ${\rG_2}$.
      \begin{enumerate}
            \item In $R(\rG_2)$, we have
            \begin{align*}
                  \widehat{J_Q(\St_{\GL_2}\nu_2^{1/2})} &= J_Q(\St_{\GL_2}\nu_2^{1/2}), \\
                  \widehat{\pi_{\deg}[1]} &= J_P(\St_{\GL_2}\nu_2^{1/2}).
            \end{align*}
            \item In $R(\PGSO_8 \rtimes \theta)$, we have
            \begin{align*}
                  \widetilde{\widehat{J_{Q'}((\St_{\GL_2} \nu_2^{1/2})_{M_{Q'}})}}
                  &= \widetilde{J_{Q'}((\St_{\GL_2} \nu_2^{1/2})_{M_{Q'}})}, \\
                  \widetilde{\widehat{\pi'_{\temp, \mathrm{ng}}}}
                  &= \widetilde{J_{P'}((\St_{\GL_2}\nu_2^{1/2})_{M_{P'}})}.
            \end{align*}
      \end{enumerate}

      We set
      \begin{align*}
            \pi_{\psi_A}(\mathbf{1}) &= J_P(\St_{\GL_2}\nu_2^{1/2}) \oplus \pi_{\mathrm{sc}}[1], \\
            \pi_{\psi_A}(\epsilon) &= J_Q(\St_{\GL_2}\nu_2^{1/2}).
      \end{align*}
\end{proposition}

\begin{proof}
      The statement $(1)$ follows from Proposition \ref{prop:subregular-abelian-Aubert-dual-computation}.
      Since the Aubert duality preserves cuspidal supports, the first equation of the statement $(2)$ follows from \cite{BanJantzen2003-DegeneratePrincipalSeriesEvenOrthogonalGroups}*{Theorem 5.3} with $\rho = \mathbf{1}, \sigma_0 = \mathbf{1}, n = 0, k = l = 2$ and $\alpha = -\frac{1}{2}$ and Proposition \ref{prop:Xu-Aubert}.
      Since by \cite{Xu2024-ArthurPacketsQuasisplitGSp2nGO2n}*{Lemma 9.2} and Lemma \ref{lem:triality-invariant-non-tempered}, the representation $J_{P'}((\St_{\GL_2}\nu_2^{1/2})_{M_{P'}})$ must be contained in the A-packet for $\PGSO_8(F)$ corresponding to $\psi'_{A}$, the second equation of the statement $(2)$ follows.
\end{proof}

\Needspace{8\baselineskip}
\begin{proposition}\label{prop:mixed-short-root-quadratic-ecr}
      We have the following results for the mixed short root parameters of ${\rG_2}$ associated with a non-trivial quadratic character $\chi$.
      \begin{enumerate}
            \item In $R(\rG_2)$, we have $\widehat{J_Q(\St_{2}(\chi)\nu_2^{1/2})} = J_P(\St_{2}(\chi)\nu_2^{1/2})$.
            \item In $R(\PGSO_8 \rtimes \theta)$, we have
            \begin{align*}
                  \widetilde{\widehat{J_{Q'}((\St_{2}(\chi) \nu_2^{1/2})_{M_{Q'}})}}
                  &= \widetilde{J_{P'}((\St_{2}(\chi) \nu_2^{1/2})_{M_{P'}})}, \\
                  \widetilde{\widehat{\pi'_{\mathrm{sc}, \chi}}}
                  &= \widetilde{\pi'_{\mathrm{sc}, \chi}}.
            \end{align*}
      \end{enumerate}

      We set
      \begin{align*}
            \pi_{\psi_A}(\mathbf{1}) &= J_P(\St_{2}(\chi)\nu_2^{1/2}), \\
            \pi_{\psi_A}(\epsilon) &= \pi_{\mathrm{sc}}[\chi].
      \end{align*}
\end{proposition}

\begin{proof}
      It suffices to show that the representation $J_{P'}((\St_{2}(\chi) \nu_2^{1/2})_{M_{P'}})$ is contained in the A-packet for $\PGSO_8(F)$ corresponding to $\psi'_{A}$.
      This follows from \cite{Xu2024-ArthurPacketsQuasisplitGSp2nGO2n}*{Lemma 9.3} and Lemma \ref{lem:triality-invariant-non-tempered}.
\end{proof}

\section{Summary}\label{sec:summary}

\begin{definition}\label{def:summary-packets}
      Let $F$ be a non-archimedean local field of characteristic zero.
      Let $\psi$ be a non-tempered A-parameter of $\rG_2(F)$.
      \begin{enumerate}
            \item We have defined a map from the set of irreducible representations of the component group $\cS_{\psi}$ to the set of finite length representations of $\rG_2(F)$ in Propositions \ref{prop:principal-ecr}, \ref{prop:subregular-trivial-ecr}, \ref{prop:subregular-quadratic-ecr}, \ref{prop:subregular-cubic-ecr}, \ref{prop:subregular-full-ecr}, \ref{prop:co-tempered-long-short-ecr}, \ref{prop:mixed-long-root-ecr}, \ref{prop:mixed-short-root-trivial-ecr}, \ref{prop:mixed-short-root-quadratic-ecr} and Proposition \ref{prop:co-tempered-long-short-unbounded-ecr}.
            For $\rho \in \Irr_{\C}(\cS_{\psi})$, let $\pi_{\psi}(\rho)$ denote the corresponding representation of $\rG_2(F)$ of finite length.
            Then, we define the A-packet $\Pi_{\psi}$ of $\rG_2(F)$ to be the set of representations
            \begin{align*}
                  \Pi_{\psi} = \bigcup_{\rho \in \Irr_{\C}(\cS_{\psi})} \mathrm{JH}(\pi_{\psi}(\rho)).
            \end{align*}
            For $s \in \cS_{\psi}$, we also set
            \begin{align*}
                  \Theta_{\psi, s} = \sum_{\rho \in \Irr_{\C}(\cS_{\psi})} \tr\rho (s s_{\psi}) \Theta_{\pi_{\psi}(\rho)}.
            \end{align*}
            \item Let $\psi'$ be the lift of $\psi$ to $\PGSO_8(F)$.
            Then, we define the A-packet $\Pi_{\psi'}$ of $\PGSO_8(F)$ to be Xu's A-packet appearing in Propositions \ref{prop:principal-ecr}, \ref{prop:subregular-trivial-ecr}, \ref{prop:subregular-quadratic-ecr}, \ref{prop:subregular-cubic-ecr}, \ref{prop:subregular-full-ecr}, \ref{prop:co-tempered-long-short-ecr}, \ref{prop:mixed-long-root-ecr}, \ref{prop:mixed-short-root-trivial-ecr}, \ref{prop:mixed-short-root-quadratic-ecr} and Proposition \ref{prop:co-tempered-long-short-unbounded-ecr}.
            Each A-packet consists of $\Aut(\PGSO_8)(F)$-invariant representations of $\PGSO_8(F)$ and there exists a map
            \begin{align*}
                  \Pi_{\psi'} \to \Irr(\cS_{\psi'}) \colon \pi' \mapsto \langle \cdot, \pi' \rangle_{\psi'},
            \end{align*}
            constructed by Arthur and Xu, see \cite{Xu2024-ArthurPacketsQuasisplitGSp2nGO2n}*{Corollary 6.2}.
      \end{enumerate}
\end{definition}

\begin{remark}
      We use similar notation for $\SO_4$ and $\PGL_3$.
      For an A-parameter $\psi$ of $\SO_4(F)$ or $\PGL_3(F)$, we denote by $\Pi_{\psi}$ the A-packet of $\SO_4(F)$ or $\PGL_3(F)$ and let $S\Theta_{\psi}$ be the stable distribution of $\SO_4(F)$ or $\PGL_3(F)$ associated with $\psi$.
\end{remark}

\Needspace{10\baselineskip}
\begin{theorem}\label{thm:main-ecr}
      Let $F$ be a non-archimedean local field of characteristic zero.
      Let $\psi$ be a non-tempered A-parameter of $\rG_2(F)$ and $\psi' = \iota \circ \psi$ be the lift of $\psi$ to $\PGSO_8(F)$.
      Then, the following statements hold.
      \begin{enumerate}
            \item We have
            \begin{align*}
                  \Pi_{\widehat{\psi}} &= \widehat{\Pi_{\psi}}, \\
                  \Pi_{\widehat{\psi'}} &= \widehat{\Pi_{\psi'}}.
            \end{align*}

            \item The distribution
            \begin{align*}
                  S\Theta_{\psi} = \sum_{\rho \in \Irr_{\C}(\cS_{\psi})} \tr\rho(s_{\psi}) \Theta_{\pi_{\psi}(\rho)}
            \end{align*}
            is stable.

            \item The twisted endoscopic character identities hold:
            \begin{align*}
                  \sum_{\pi' \in \Pi_{\psi'}} \langle s s_{\psi'}, \pi' \rangle_{\psi'} \Theta_{\widetilde{\pi'}}
                  = \Trans^{\widetilde{\PGSO_8}, *}_{H_s} S\Theta^{H_s}_{\psi_s},
            \end{align*}
            where $s \in \cS_{\psi'}$, $H_s$ is the endoscopic group of $\PGSO_8 \rtimes \theta$ associated with $s \rtimes \theta$, $\psi_s$ is the A-parameter of $H_s$ associated with $\psi'$, and $S\Theta^{H_s}_{\psi_s}$ is the stable distribution of $H_s$ associated with $\psi_s$.

            \item The endoscopic character identities for $A$-parameters $\psi$ of $\rG_2(F)$ hold:

            \begin{align*}
                  \Theta_{\psi, s}
                  = \Trans^{\rG_2, *}_{H_s} S\Theta^{H_s}_{\psi_s},
            \end{align*}
            where $s \in \cS_{\psi}$, $H_s$ is the endoscopic group of $\rG_2(F)$ associated with $s$, $\psi_s$ is the A-parameter of $H_s$ associated with $\psi$, and $S\Theta^{H_s}_{\psi_s}$ is the stable distribution of $H_s$ associated with $\psi_s$.
      \end{enumerate}
\end{theorem}

\begin{proof}
      These follow from Propositions \ref{prop:principal-ecr}, \ref{prop:subregular-trivial-ecr}, \ref{prop:subregular-quadratic-ecr}, \ref{prop:subregular-cubic-ecr}, \ref{prop:subregular-full-ecr}, \ref{prop:co-tempered-long-short-ecr}, \ref{prop:mixed-long-root-ecr}, \ref{prop:mixed-short-root-trivial-ecr}, \ref{prop:mixed-short-root-quadratic-ecr} and Proposition \ref{prop:co-tempered-long-short-unbounded-ecr}.
\end{proof}

\section*{Acknowledgements}
The author thanks Sug Woo Shin for his suggestion to follow the approach of Atobe--Gan--Ichino--Kaletha--M\'inguez--Shin \cite{AtobeGanEtAl2025-LocalIntertwiningRelationsCoTemperedPackets} and for his careful comments on the draft.

\section*{AI Declaration}
The author wrote the initial draft without using AI and subsequently used ChatGPT (OpenAI, Astra 6) to revise the language and organization, review references and mathematical exposition, prepare the proof diagram, and edit the LaTeX source.
The author takes full responsibility for the content of this paper.

\bibliographystyle{plain}
\bibliography{bibliography}
\end{document}

%% file: proof-of-main-theorem.tex
\hypersetup{linkcolor=black}
\newcommand{\proofnode}[2]{{\sffamily\bfseries #1}\\[3pt]#2}
\begin{tikzpicture}[
    x=0.9cm,
    every node/.style={text=black},
    box/.style={draw=black, rounded corners=2pt, align=center,
        text width=29mm, inner xsep=2.2mm, inner ysep=2.4mm,
        minimum height=17mm, font=\fontsize{10}{12}\selectfont},
    family/.style={box, text width=39.5mm, fill=white},
    wire/.style={draw=black, line width=0.65pt, line cap=round},
    arrow/.style={wire, -{Stealth[length=2mm,width=1.5mm]}},
    junction dot/.style={circle, fill=black, inner sep=0pt, minimum size=2.5pt}
]
\node[box] (xu) at (-6.15,-0.5)
    {\proofnode{Proposition~\ref{prop:Xu-Aubert}}{Aubert duality for\\Xu packets}};
\node[box, text width=32mm] (temp) at (-2.05,-0.5)
    {\proofnode{Theorem~\ref{thm:tempered-input}}{Tempered endoscopic\\character identities}};
\node[box] (sign) at (2.05,-0.5)
    {\proofnode{Theorem~\ref{thm:signed-transfer}}{Triality transfer\\and Aubert signs}};
\node[box, text width=32mm] (lift) at (6.15,-0.5)
    {\proofnode{Lemma~\ref{lem:mixed-long-lift}\\Corollary~\ref{cor:mixed-long-exact}}{Mixed long-root lift\\and exact sequences}};

\coordinate (shared) at (0,-2.15);
\coordinate (fork) at (0,-2.45);
\coordinate (cj) at (-5.35,-3.35);
\coordinate (mj) at (0.25,-3.35);
\draw[wire] (temp.south) .. controls +(0,-0.53) and +(0,0.42) .. (shared);
\draw[wire] (sign.south) .. controls +(0,-0.53) and +(0,0.42) .. (shared);
\draw[wire] (shared) -- (fork);
\draw[wire] (fork) .. controls +(0,-0.48) and +(0,0.56) .. (cj);
\draw[wire] (fork) .. controls +(0,-0.38) and +(0,0.48) .. (mj);
\draw[wire] (xu.south) .. controls +(0,-0.93) and +(0,0.76) .. (cj);
\draw[wire] (lift.south) .. controls +(0,-1.22) and +(0,0.83) .. (mj);
\foreach \p in {shared,fork,cj,mj} \node[junction dot] at (\p) {};

\node[family] (co) at (-5.35,-4.75)
    {\proofnode{Corollary~\ref{cor:cotempered}}{Co-tempered packets\\from tempered packets}};
\node[family] (long) at (0.25,-4.75)
    {\proofnode{Proposition~\ref{prop:mixed-long-root-ecr}}{Mixed long-root packets\\and character identities}};
\node[box, text width=37mm] (levi) at (5.8,-4.75)
    {\proofnode{Lemma~\ref{lem:induction-transfer} and Mui\'c}{Transfer under\\induction and irreducibility}};
\draw[arrow] (cj) -- (co.north);
\draw[arrow] (mj) -- (long.north);

\node[family] (cases) at (-5.35,-7.55)
    {\proofnode{Sections~\ref{sec:principal-parameter}--\ref{sec:co-tempered-long-short-supercuspidal}}{Principal and subregular cases\\Root cases attached to\\supercuspidals of \(\mathrm{PGL}_2(F)\)}};
\node[family] (short) at (0.25,-7.55)
    {\proofnode{Propositions~\ref{prop:mixed-short-root-trivial-ecr}--\ref{prop:mixed-short-root-quadratic-ecr}}{Mixed short-root packets}};
\node[family, text width=37mm] (series) at (5.8,-7.55)
    {\proofnode{Proposition~\ref{prop:co-tempered-long-short-unbounded-ecr}}{Root packets from\\unitary principal series\\of \(\mathrm{PGL}_2(F)\)}};
\draw[arrow] (co.south) -- (cases.north);
\draw[arrow] (long.south) -- node[right,font=\fontsize{9}{11}\selectfont,inner sep=3pt]
    {Aubert duality} (short.north);
\draw[arrow] (levi.south) -- (series.north);

\coordinate (collect) at (0.25,-9.6);
\draw[wire] (cases.south) .. controls +(0,-0.91) and +(0,0.62) .. (collect);
\draw[wire] (series.south) .. controls +(0,-0.91) and +(0,0.62) .. (collect);
\coordinate (mixed) at (0.25,-8.92);
\draw[wire] (long.west) .. controls +(-0.49,0) and +(0,0.64) .. (-2.53,-6.27)
    -- (-2.53,-8.43) .. controls +(0,-0.43) and +(0,0.32) .. (mixed);
\draw[wire] (short.south) -- (mixed) -- (collect);
\node[junction dot] at (mixed) {};
\node[junction dot] at (collect) {};
\node[box, fill=white, line width=1pt, draw=black, text width=112mm] (main) at (0.25,-11.05)
    {\proofnode{Definition~\ref{def:summary-packets} and Theorem~\ref{thm:main-ecr}}
      {Arthur packets, compatibility with Aubert duality, stability,\\
       and the ordinary and twisted endoscopic character identities}};
\draw[arrow] (collect) -- (main.north);
\end{tikzpicture}

%% file: bibliography.bib
@unpublished{AtobeGanEtAl2025-LocalIntertwiningRelationsCoTemperedPackets,
  AUTHOR = {Atobe, Hiraku and Gan, Wee Teck and Ichino, Atsushi and
            Kaletha, Tasho and M{\'i}nguez, Alberto and Shin, Sug Woo},
  TITLE = {Local intertwining relations and co-tempered {$A$}-packets
           of classical groups},
  YEAR = {2025},
  NOTE = {Preprint, \href{https://arxiv.org/abs/2410.13504v2}{arXiv:2410.13504v2}},
  eprint = {2410.13504v2},
  archivePrefix = {arXiv},
  primaryClass = {math.NT},
  url = {https://arxiv.org/abs/2410.13504v2},
}

@article {BanJantzen2003-DegeneratePrincipalSeriesEvenOrthogonalGroups,
    AUTHOR = {Ban, Dubravka and Jantzen, Chris},
     TITLE = {Degenerate principal series for even-orthogonal groups},
   JOURNAL = {Represent. Theory},
  FJOURNAL = {Representation Theory. An Electronic Journal of the American
              Mathematical Society},
    VOLUME = {7},
      YEAR = {2003},
     PAGES = {440--480},
      ISSN = {1088-4165},
   MRCLASS = {22E50},
  MRNUMBER = {2017065},
MRREVIEWER = {David\ Manderscheid},
       DOI = {10.1090/S1088-4165-03-00166-3},
       URL = {https://doi-org.kyoto-u.idm.oclc.org/10.1090/S1088-4165-03-00166-3},
}

@unpublished{ChengKaletha2026-AubertBernsteinInvolutionsDisconnectedGroups,
  author = {Cheng, Yongshen and Kaletha, Tasho},
  title = {The {A}ubert and {B}ernstein involutions for disconnected groups},
  year = {2026},
  note = {Preprint, \href{https://arxiv.org/abs/2607.10660}{arXiv:2607.10660}},
  eprint = {2607.10660},
  archivePrefix = {arXiv},
  primaryClass = {math.RT},
  url = {https://arxiv.org/abs/2607.10660}
}

@article{GanGurevichJiang2002-CubicUnipotentArthurParameters,
  AUTHOR = {Gan, Wee Teck and Gurevich, Nadya and Jiang, Dihua},
  TITLE = {Cubic unipotent {A}rthur parameters and multiplicities of
           square integrable automorphic forms},
  JOURNAL = {Invent. Math.},
  VOLUME = {149},
  NUMBER = {2},
  YEAR = {2002},
  PAGES = {225--265},
  MRNUMBER = {1918673},
  DOI = {10.1007/s002220200210},
  URL = {https://doi.org/10.1007/s002220200210}
}

@article {GanSavin2023-HoweDualityDichotomyExceptionalThetaCorrespondences,
    AUTHOR = {Gan, Wee Teck and Savin, Gordan},
     TITLE = {Howe duality and dichotomy for exceptional theta
              correspondences},
   JOURNAL = {Invent. Math.},
  FJOURNAL = {Inventiones Mathematicae},
    VOLUME = {232},
      YEAR = {2023},
    NUMBER = {1},
     PAGES = {1--78},
      ISSN = {0020-9910,1432-1297},
   MRCLASS = {11F27 (11F70 22E50)},
  MRNUMBER = {4557399},
       DOI = {10.1007/s00222-022-01165-2},
       URL = {https://doi.org/10.1007/s00222-022-01165-2},
}

@article {GanSavin2023-LocalLanglandsConjectureG2,
    AUTHOR = {Gan, Wee Teck and Savin, Gordan},
     TITLE = {The local {L}anglands conjecture for {$G_2$}},
   JOURNAL = {Forum Math. Pi},
  FJOURNAL = {Forum of Mathematics. Pi},
    VOLUME = {11},
      YEAR = {2023},
     PAGES = {Paper No. e28, 42},
      ISSN = {2050-5086},
   MRCLASS = {11S37 (11F27 11F70 22E50)},
  MRNUMBER = {4658199},
MRREVIEWER = {Yiwen\ Ding},
       DOI = {10.1017/fmp.2023.27},
       URL = {https://doi.org/10.1017/fmp.2023.27},
}

@unpublished{Kaletha2022-LocalLanglandsConjecturesDisconnectedGroups,
  author = {Kaletha, Tasho},
  title = {On the local {L}anglands conjectures for disconnected groups},
  year = {2022},
  note = {Preprint, version 1, \href{https://arxiv.org/abs/2210.02519v1}{arXiv:2210.02519}},
  eprint = {2210.02519v1},
  archivePrefix = {arXiv},
  primaryClass = {math.RT},
  url = {https://arxiv.org/abs/2210.02519v1}
}

@book {MoeglinWaldspurger2016-StabilisationFormuleTracesTordue1,
    AUTHOR = {M{\oe}glin, Colette and Waldspurger, Jean-Loup},
     TITLE = {Stabilisation de la formule des traces tordue. {V}ol. 1},
    SERIES = {Progress in Mathematics},
    VOLUME = {316},
 PUBLISHER = {Birkh\"auser/Springer, Cham},
      YEAR = {2016},
     PAGES = {xviii+587},
      ISBN = {978-3-319-30048-1; 978-3-319-30049-8},
   MRCLASS = {22E50 (11F70 11F72 22E55)},
  MRNUMBER = {3823813},
MRREVIEWER = {Jean-Pierre\ Labesse},
       DOI = {10.1007/978-3-319-30049-8},
       URL = {https://doi.org/10.1007/978-3-319-30049-8},
}

@article{Qin2026-AubertZelevinskyInvolutionG2HeckeAlgebras,
  author = {Qin, Chuan},
  title = {The {A}ubert-{Z}elevinsky involution for {$G_2$} and its associated {H}ecke algebras},
  journal = {J. Pure Appl. Algebra},
  volume = {230},
  number = {11},
  year = {2026},
  pages = {108372},
  MRNUMBER = {5118070},
  doi = {10.1016/j.jpaa.2026.108372},
  url = {https://doi.org/10.1016/j.jpaa.2026.108372}
}

@unpublished{TakanashiEndoscopyLocal,
  author = {Takanashi, Yugo},
  title = {Endoscopic description of the local {L}anglands correspondence for {$G_2$}},
  year = {2026},
  note = {Preprint, \href{https://arxiv.org/abs/2609.19439v1}{arXiv:2609.19439v1}}
}

@unpublished{TakanashiThesisLocal,
  author = {Takanashi, Yugo},
  title = {On the {H}iraga--{I}chino--{I}keda conjecture on formal degrees for {$G_2$}},
  year = {2026},
  note = {Preprint, \href{https://arxiv.org/abs/2609.19430v1}{arXiv:2609.19430v1}}
}

@article{Xu2017-MoeglinParametrizationArthurPackets,
  author = {Xu, Bin},
  title = {On {M}\oe glin's parametrization of {A}rthur packets for {$p$}-adic quasisplit {$\mathrm{Sp}(N)$} and {$\mathrm{SO}(N)$}},
  journal = {Canad. J. Math.},
  volume = {69},
  number = {4},
  year = {2017},
  pages = {890--960},
  MRNUMBER = {3679701},
  doi = {10.4153/CJM-2016-029-3},
  url = {https://doi.org/10.4153/CJM-2016-029-3}
}

@article{Xu2017-CuspidalSupportDiscreteSeries,
  author = {Xu, Bin},
  title = {On the cuspidal support of discrete series for {$p$}-adic
           quasisplit {$\mathrm{Sp}(N)$} and {$\mathrm{SO}(N)$}},
  journal = {Manuscripta Math.},
  volume = {154},
  number = {3-4},
  year = {2017},
  pages = {441--502},
  MRNUMBER = {3713922},
  doi = {10.1007/s00229-017-0923-x},
  url = {https://doi.org/10.1007/s00229-017-0923-x}
}

@article {Xu2018-LPacketsQuasisplitGSp2nGO2n,
    AUTHOR = {Xu, Bin},
     TITLE = {L-packets of quasisplit {$\mathrm{GSp}(2n)$} and {$\mathrm{GO}(2n)$}},
   JOURNAL = {Math. Ann.},
  FJOURNAL = {Mathematische Annalen},
    VOLUME = {370},
      YEAR = {2018},
    NUMBER = {1-2},
     PAGES = {71--189},
      ISSN = {0025-5831,1432-1807},
   MRCLASS = {22E50 (11F70)},
  MRNUMBER = {3747484},
MRREVIEWER = {Ivan\ Mati\'{c}},
       DOI = {10.1007/s00208-016-1515-x},
       URL = {https://doi.org/10.1007/s00208-016-1515-x},
}

@article {Xu2024-ArthurPacketsQuasisplitGSp2nGO2n,
    AUTHOR = {Xu, Bin},
     TITLE = {Arthur packets for quasisplit {${\rm GSp}(2n)$} and {${\rm GO}(2n)$} over
              a {$p$}-adic field},
   JOURNAL = {Int. Math. Res. Not. IMRN},
  FJOURNAL = {International Mathematics Research Notices. IMRN},
      YEAR = {2024},
    NUMBER = {13},
     PAGES = {10349--10409},
      ISSN = {1073-7928,1687-0247},
   MRCLASS = {22E50 (11F70 11S37)},
  MRNUMBER = {4770373},
MRREVIEWER = {Amiya\ Kumar\ Mondal},
       DOI = {10.1093/imrn/rnae086},
       URL = {https://doi-org.kyoto-u.idm.oclc.org/10.1093/imrn/rnae086},
}

@article{HuangMagaardSavin1998-UnipotentRepresentationsG2,
  author = {Huang, Jing-Song and Magaard, Kay and Savin, Gordan},
  title = {Unipotent representations of {$G_2$} arising from the minimal
           representation of {$D_4^E$}},
  journal = {J. Reine Angew. Math.},
  volume = {500},
  year = {1998},
  pages = {65--81},
  MRNUMBER = {1637485},
  doi = {10.1515/crll.1998.071},
  url = {https://doi.org/10.1515/crll.1998.071}
}

@article {Muic1997-UnitaryDualG2,
    AUTHOR = {Mui\'{c}, Goran},
     TITLE = {The unitary dual of {$p$}-adic {$G_2$}},
   JOURNAL = {Duke Math. J.},
  FJOURNAL = {Duke Mathematical Journal},
    VOLUME = {90},
      YEAR = {1997},
    NUMBER = {3},
     PAGES = {465--493},
      ISSN = {0012-7094,1547-7398},
   MRCLASS = {22E50},
  MRNUMBER = {1480543},
MRREVIEWER = {David\ Manderscheid},
       DOI = {10.1215/S0012-7094-97-09012-8},
       URL = {https://doi.org/10.1215/S0012-7094-97-09012-8},
}

@article {Hiraga2004-FunctorialityZelevinskiInvolutions,
    AUTHOR = {Hiraga, Kaoru},
     TITLE = {On functoriality of {Z}elevinski involutions},
   JOURNAL = {Compos. Math.},
  FJOURNAL = {Compositio Mathematica},
    VOLUME = {140},
      YEAR = {2004},
    NUMBER = {6},
     PAGES = {1625--1656},
      ISSN = {0010-437X,1570-5846},
   MRCLASS = {22E50 (11F70 22E35)},
  MRNUMBER = {2098406},
MRREVIEWER = {Ernst-Wilhelm\ Zink},
       DOI = {10.1112/S0010437X04000892},
       URL = {https://doi.org/10.1112/S0010437X04000892},
}

@article {SchneiderStuhler1997-RepresentationTheorySheaves,
    AUTHOR = {Schneider, Peter and Stuhler, Ulrich},
     TITLE = {Representation theory and sheaves on the {B}ruhat-{T}its
              building},
   JOURNAL = {Inst. Hautes \'{E}tudes Sci. Publ. Math.},
  FJOURNAL = {Institut des Hautes \'{E}tudes Scientifiques. Publications
              Math\'{e}matiques},
    NUMBER = {85},
      YEAR = {1997},
     PAGES = {97--191},
      ISSN = {0073-8301,1618-1913},
   MRCLASS = {22E50 (11F70 20G25)},
  MRNUMBER = {1471867},
MRREVIEWER = {Ernst-Wilhelm\ Zink},
       DOI = {10.1007/BF02699536},
       URL = {https://www.numdam.org/item/PMIHES_1997__85__97_0/},
}

@article {GanGurevich2006-NontemperedAPacketsG2,
    AUTHOR = {Gan, Wee Teck and Gurevich, Nadya},
     TITLE = {Nontempered {A}-packets of {$G_2$}: liftings from
              {$\widetilde{\rm SL}_2$}},
   JOURNAL = {Amer. J. Math.},
  FJOURNAL = {American Journal of Mathematics},
    VOLUME = {128},
      YEAR = {2006},
    NUMBER = {5},
     PAGES = {1105--1185},
      ISSN = {0002-9327,1080-6377},
   MRCLASS = {22E55 (11F70)},
  MRNUMBER = {2262172},
MRREVIEWER = {Dihua\ Jiang},
       DOI = {10.1353/ajm.2006.0040},
       URL = {https://doi.org/10.1353/ajm.2006.0040},
}

@unpublished{KottwitzShelstad2012-SplittingInvariantsSignConventions,
  author = {Kottwitz, Robert E. and Shelstad, Diana},
  title = {On splitting invariants and sign conventions in endoscopic transfer},
  year = {2012},
  note = {Preprint, \href{https://arxiv.org/abs/1201.5658}{arXiv:1201.5658}},
  url = {https://arxiv.org/abs/1201.5658}
}

@article {Muic2005-ReducibilityStandardRepresentations,
    AUTHOR = {Mui\'c, Goran},
     TITLE = {Reducibility of standard representations},
   JOURNAL = {Pacific J. Math.},
  FJOURNAL = {Pacific Journal of Mathematics},
    VOLUME = {222},
      YEAR = {2005},
    NUMBER = {1},
     PAGES = {133--168},
      ISSN = {0030-8730,1945-5844},
   MRCLASS = {22E50 (11F70 22E35)},
  MRNUMBER = {2200248},
MRREVIEWER = {Dubravka\ Ban},
       DOI = {10.2140/pjm.2005.222.133},
       URL = {https://doi-org.kyoto-u.idm.oclc.org/10.2140/pjm.2005.222.133},
}

@article {Muic2005-ReducibilityGeneralizedPrincipalSeries,
    AUTHOR = {Mui\'c, Goran},
     TITLE = {Reducibility of generalized principal series},
   JOURNAL = {Canad. J. Math.},
  FJOURNAL = {Canadian Journal of Mathematics. Journal Canadien de
              Math\'ematiques},
    VOLUME = {57},
      YEAR = {2005},
    NUMBER = {3},
     PAGES = {616--647},
      ISSN = {0008-414X,1496-4279},
   MRCLASS = {22E50 (11F70 22E35)},
  MRNUMBER = {2134404},
MRREVIEWER = {David\ Manderscheid},
       DOI = {10.4153/CJM-2005-025-4},
       URL = {https://doi-org.kyoto-u.idm.oclc.org/10.4153/CJM-2005-025-4},
}

@article {Xu2021-NonarchimedeanComponentsNonEndoscopicAutomorphicRepresentations,
    AUTHOR = {Xu, Bin},
     TITLE = {Nonarchimedean components of non-endoscopic automorphic
              representations for quasisplit {$Sp(N)$} and {$O(N)$}},
   JOURNAL = {Math. Z.},
  FJOURNAL = {Mathematische Zeitschrift},
    VOLUME = {297},
      YEAR = {2021},
    NUMBER = {1-2},
     PAGES = {885--921},
      ISSN = {0025-5874,1432-1823},
   MRCLASS = {22E50 (11F70)},
  MRNUMBER = {4204718},
MRREVIEWER = {Ivan\ Mati\'c},
       DOI = {10.1007/s00209-020-02539-z},
       URL = {https://doi-org.kyoto-u.idm.oclc.org/10.1007/s00209-020-02539-z},
}

@unpublished{HalleckDube2026-WeightedFundamentalLemmaNonSplitGroups,
  author = {Halleck-Dub\'e, Cass},
  title = {The weighted fundamental lemma for non-split groups},
  year = {2026},
  note = {Preprint, \href{https://arxiv.org/abs/2611.26750v1}{arXiv:2611.26750v1}},
  eprint = {2611.26750v1},
  archivePrefix = {arXiv},
  primaryClass = {math.RT},
  url = {https://arxiv.org/abs/2611.26750v1},
}

@book{Arthur2013-EndoscopicClassificationRepresentations,
  author = {Arthur, James},
  title = {The endoscopic classification of representations: Orthogonal and symplectic groups},
  series = {American Mathematical Society Colloquium Publications},
  volume = {61},
  publisher = {American Mathematical Society, Providence, RI},
  year = {2013},
  pages = {xviii+590},
  isbn = {978-0-8218-4990-3},
  MRNUMBER = {3135650},
  doi = {10.1090/coll/061},
  url = {https://doi.org/10.1090/coll/061}
}

@unpublished{CunninghamEtAl2026-GeometrizationKottwitzShelstadTransfer,
  AUTHOR = {Cunningham, Clifton and Hazeltine, Alexander and Liu, Baiying and
            Lo, Chi-Heng and Ray, Mishty and Xu, Bin and Sch{\"u}rmann, J{\"o}rg},
  TITLE = {Geometrization of {K}ottwitz--{S}helstad transfer and
           {V}ogan's conjecture on local {A}rthur packets},
  YEAR = {2026},
  NOTE = {Preprint, \href{https://arxiv.org/abs/2611.03342v1}{arXiv:2611.03342v1}},
  eprint = {2611.03342v1},
  archivePrefix = {arXiv},
  primaryClass = {math.RT},
  URL = {https://arxiv.org/abs/2611.03342v1},
}
